\documentclass[3p]{elsarticle}
\usepackage{dsfont}
\usepackage{amssymb}
\usepackage{amsthm}
\usepackage{amsmath}
\usepackage{lineno}
\usepackage{empheq}
\usepackage{graphicx}
\usepackage{booktabs}
\usepackage{lineno}

\newdefinition{rmk}{Remark}
\newproof{pf}{Proof}
\newproof{pot}{Proof of Theorem \ref{thm2}}
\usepackage{caption}
\usepackage{subcaption}
\begin{document}
\begin{frontmatter}
\title{The Laguerre finite difference one-way equation solver}
\author{Andrew V. Terekhov}
\ead{andrew.terekhov@mail.ru}
\address{Institute of
Computational Mathematics and Mathematical Geophysics,
630090, Novosibirsk, Russia}
\address{Novosibirsk State University, 630090, Novosibirsk,
Russia}

\begin{abstract}
This paper presents a new finite difference algorithm for solving the 2D one-way wave equation with a preliminary approximation of a pseudo-differential operator by a system of partial differential equations. As opposed to the existing approaches, the integral Laguerre  transform  instead of Fourier transform is used. After carrying out the approximation of spatial variables it is possible to obtain systems of linear algebraic equations with better computing properties and to reduce computer costs for their solution.  High accuracy of calculations is attained at the expense of employing finite difference approximations of higher accuracy order that are based on the dispersion-relationship-preserving method and the Richardson extrapolation in the downward  continuation direction. The numerical experiments  have verified  that as compared to the spectral difference method based on Fourier transform, the new algorithm allows one to calculate wave fields with a higher  degree of accuracy and a lower level of numerical noise and artifacts including those for non-smooth velocity models. In the context of solving the geophysical problem the post-stack migration for velocity models of the types Syncline and Sigsbee2A has been carried out. It is shown that the images obtained contain lesser noise and are considerably better focused as compared to those obtained by the known Fourier Finite Difference and Phase-Shift Plus Interpolation methods. There is an opinion that purely finite difference approaches do not allow carrying out the seismic migration procedure with sufficient accuracy, however the results obtained disprove this statement. For the supercomputer implementation it is proposed to use the parallel dichotomy algorithm when solving systems of linear algebraic equations with block-tridiagonal matrices.
\end{abstract}
\end{frontmatter}
\section{Introduction}
The wave equation describes waves propagating in all directions.
A one-way wave equation, also known as  paraxial or parabolic, wave equation is an equation describing only downward propagating waves or only upward
propagating waves for a positive or a negative vertical component of propagation velocity, respectively.
Mathematical models based on the wave equation are often considered in problems of seismic prospecting \cite{Claerbout1971,seis_interpr,Angus2013}, ocean acoustics \cite{Tappert1977,Jensen2011,Lee2000} as well as for setting non-reflecting boundary conditions \cite{Lindman1975,Engquist1977,Trefethen1986}, whereas initially the ideas were formulated for solving the problem of electromagnetic wave propagation along the Earth's surface\cite{Leontovich1946}.

For solving the one-way equation it is possible to use finite difference methods \cite{Claerbout:1985,Lee1982}. To this end, the pseudo-differential square-root operator was first approximated by an optimized series which has its origin in a continued fraction expansion \cite{Claerbout:1985,Halpern1988,Bamberger1988,Bamberger1988a}.  A finite difference method can handle strong lateral velocity variations, however it is dip-limited and has significant numerical dispersions and artifacts. A conventional approach to reducing numerical dispersions is in applying finite difference schemes of a high order of accuracy, but when solving the one-way equation there arise essential computational problems. First, there is a numerical instability of calculations in using schemes of high orders of accuracy for the approximation of spatial derivatives towards the wave field extrapolation. Another difficulty is in that for decreasing the time of calculation, a practical approach, independently proposed in \cite{Marchuk1968} and in \cite{Strang1968}, to solving extrapolation equation (\ref{pade_approx}) is implemented by the Marchuk-Strang splitting \cite{Marchuk1990}. However for an inhomogeneous velocity model such an approach provides only the first accuracy order as related to a step in depth because the operators to be split are not mutually commutative.
In order to attain second order of accuracy, the splitting procedure should be specially symmetrized  \cite{Marchuk1990} thus twice increasing computer costs. The abandonment of using the Marchuk-Strang splitting procedure as regards this algorithm demands higher computer costs, especially for a 3D case when it is necessary to use iterative procedures for solving systems of linear algebraic equations (SLAEs). Another approach to reducing numerical dispersion and computing artifacts when applying the finite difference approximation is in the use of the filtration procedures \cite{Bunks1995} as well as of the nonlinear procedures \cite{Fei1995} similar to the TVD-schemes concepts for hyperbolic equations.

For reducing the numerical dispersion when solving geophysics problems the Fourier-based methods were developed \cite{Gazdag1978,Gazdag1984,Ristow01121994}.  These methods only approximately take into account the lateral velocity variations as being obtained on the basis of non-identical mathematical transformations thus unpredictedly affecting the quality of solution. Nevertheless, for real practical tasks it has been just these methods that make possible to carry out calculations with sufficiently large mesh sizes. In this case numerical artifacts generated by the Fourier-based methods are essentially smaller than those generated by finite difference methods of second order of accuracy. On the other hand, the advantage of the finite difference approximation is in that it possesses convergence and is spatially localized, thus allowing the evaluation of the quality of solutions that are obtained with a sequence of imbedded meshes. In \cite{Guddati2006}, an alternative approach to obtaining solutions similar to that of the one-way equation is proposed. This approach is a combination of various ideas related to the one-way equations, half-space stiffness relation, special finite element discretization and complex coordinate stretching. However, the question of constructing a stable difference method of a high accuracy order for the classical one-way equation is still an open question and urgent both from the standpoint of applications and methods.

In this paper we propose a new finite difference algorithm of high order of accuracy for solving the one-way equation. As opposed to the generally accepted approach we use the integral Laguerre transform \cite{abramowitz+stegun} instead of Fourier transform with respect to time. As for the approximation of spatial derivatives we employ only finite difference schemes. The approach in question was investigated in solving dynamic problems for the acoustics and elasticity equations \cite{Mikh2003,fatab2011,Terekhov:2013,Terekhov2015206}. Applying the Laguerre transform as well as Fourier transform with respect to time enables us not to consider  time-centering of the method of calculating the wave field that is required in analyzing finite difference approximations of the temporal derivative. Also, the Laguerre transform makes possible to obtain SLAEs with considerably better computational properties, than with using Fourier transform. Decreasing computer costs for solving SLAEs allows the abandonment of the Marchuk-Strang splitting and hence an increase both in stability and accuracy of the algorithm as a whole. In this case the approximation error relatively the mesh size towards the extrapolation of the wave field is of fourth accuracy order thus considerably reducing the numerical dispersion as compared to algorithms accurate to second order. As a result, a widespread opinion that a purely finite difference approach does not allow obtaining high quality solutions to the one-way equation will be invalidated.
\section{Description of the method}
\subsection{Governing  wave equations}
The one-way equation can be obtained from the acoustics equation
\begin{equation}
\label{waveeq}
\frac{1}{c^2}\frac{\partial^2 u}{\partial t^2 }-\nabla^2u=0,
\end{equation}
where $u\equiv u(x,z,t)$ is the field variable, $c$ is the wave velocity, the vertical direction $z$ is the extrapolation direction, i.e., the direction of one-way propagation and the positive axis $z$ is directed downward, i.e., toward increasing depth.
Applying Fourier transforms to equation (\ref{waveeq}) over $x$ and $t$ yields
\begin{equation}
\frac{\partial^2 \tilde{u}}{\partial z^2 }+\left(\frac{\omega^2}{c^2}-k_x^2\right)\tilde{u}=0,
\end{equation}
where $\tilde{u}\equiv \tilde{u}(x,z,\omega)$ is a wave component at the radial frequency $\omega$, $k_x$ is the horizontal wave number. Then, the 2D one-way acoustic wave equation for upcoming waves in the frequency
domain can be expressed as:
\begin{equation}
\frac{\partial \tilde{u}}{\partial z}=-\mathrm{i}\frac{\omega}{c}\sqrt{1-\left(\frac{c k_x}{\omega}\right)^2} \tilde{u},
\label{one_way}
\end{equation}
where $\mathrm{i}=\sqrt{-1}$ is the imaginary unit.
In practice, the square-root operator is approximated by an optimized series which has its origin in a continued fraction expansion that can be represented by the ratios of polynomials \cite{Claerbout:1985}.
\begin{equation}
\label{pade_approx}
\frac{\partial \tilde{u}}{\partial z}=-\mathrm{i}\frac{\omega}{c}\left[1-\sum_{s=1}^n\frac{\beta_sc^2 k_x^2}{\omega^2-\gamma_sc^2 k_x^2}\right]\tilde{u},
\end{equation}
and the polynomial coefficients $\gamma_s,\beta_s \in \mathfrak{R}$ for the propagation angle should be optimized\cite{Halpern1988,Lee01101985}.

Let us write down the one-way equation for the spatial-temporal domain $(x,z,t)$. To this end it is convenient to introduce the auxiliary functions\cite{Bamberger1988,Bamberger1988a} of the form  $$\tilde{\psi}_s=\frac{\beta_sc^2k_x^2}{\omega^2-\gamma_s c^2k_x^2}\tilde{u},$$
so that equation (\ref{pade_approx}) can be rewritten as
$$
c\frac{\partial \tilde{u}}{\partial z}+i\omega \tilde{u}-i\omega\sum_{s=1}^{n}\tilde{\psi}_s=0.
$$
Applying the inverse Fourier transform to the two latter equations, we come to a system of equations determining the down-continuation process
\begin{subequations}
\label{main_eq}
\begin{empheq}[left=\empheqlbrace]{align}
\label{main_eq1}
\frac{\partial u}{\partial t}+c\frac{\partial u}{\partial z}-\sum_{s=1}^{n}\frac{\partial \psi_s}{\partial t}=0,\\
\frac{1}{c^2}\frac{\partial^2 \psi_s}{\partial t^2}-\gamma_s \frac{\partial^2\psi_s }{\partial x^2}-\beta_s\frac{\partial^2 u }{\partial x^2}=0,\quad s=1,2,...n.
\label{main_eq2}
\end{empheq}
\end{subequations}

\subsection{Temporal approximation}
Let us consider the direct and the inverse Laguerre transforms\cite{Integral_Transform}
\begin{equation}
\label{series_lag}
L\{g(t)\}=\bar{g}_m=\int_{0}^{\infty}\left(\eta t\right)^{-\alpha/2}g(t)l^{\alpha}_m(\eta t)dt,\quad\quad g(t)=L^{-1}\{\bar{g}_m\}=(\eta t)^{\frac{\alpha}{2}}\sum_{m=0}^{\infty}\bar{g}_ml^{\alpha}_m(\eta t).
\end{equation}
Here $l^{\alpha}_m(t)$ are the orthonormal Laguerre functions, $m$ is Laguerre polynomial degree, $\alpha$ is the order of the Laguerre functions, and $\eta$ is the transformation parameter. Assuming \mbox{$\lim_{t\rightarrow \infty}g(t)=\lim_{t\rightarrow \infty}\frac{dg}{dt} (t)=0$}, we can show \cite{Integral_Transform,Mikhailenko1999} that
\begin{equation}
L\left\{\frac{d^{k}}{d t^k}g(t)\right\}=\left(\frac{\eta}{2}\right)^{k}\bar{g}_m+\Phi_k(\bar{g}_m),
\label{partial_t_lag}
\end{equation}
where for $k=1,2$ we have
\begin{equation}
\label{lag33}
\Phi_1(\bar{g}_n)\equiv \eta\sqrt{\frac{n!}{(n+\alpha)!}}\sum_{k=0}^{n-1}\sqrt{\frac{(k+\alpha)!}{k!}}\bar{g}_k, \quad\Phi_2(\bar{g}_n)\equiv\eta^2\sqrt{\frac{n!}{(n+\alpha)!}}\sum_{k=0}^{n-1}(n-k)\sqrt{\frac{(k+\alpha)!}{k!}}\bar{g}_k.
\end{equation}
Then applying the Laguerre transform (\ref{series_lag})  to equations (\ref{main_eq1}),(\ref{main_eq2})  we obtain the following system of equations for the calculation of the $m$-th coefficient of expansion:
\begin{subequations}
\label{main_eq_laguerre}
\begin{empheq}[left=\empheqlbrace]{align}
\label{main_eq1_laguerre}
\tilde{\eta} \bar{u}^m+c\frac{\partial \bar{u}^m}{\partial z}=\sum_{s=1}^3\left(\tilde{\eta}\bar{\psi}_s^m+\Phi_1\left(\bar{\psi}_s^m\right)\right)-\Phi_1\left(\bar{u}^m\right),\\
\label{main_eq2_laguerre}
c^2\gamma_s\frac{\partial^2\bar{\psi}_s^m }{\partial x^2}-\tilde{\eta}^2\bar{\psi}_s^m+\beta_sc^2\frac{\partial^2\bar{u}^m }{\partial x^2}=\Phi_2(\bar{\psi}_s^m) ,\quad s=1,2,3,
\end{empheq}
\end{subequations}
where $\tilde{\eta}=\eta/2$, and the index $m$ denotes number of a term of series (\ref{series_lag}).

The expansion parameters $\gamma,\beta$ from (\ref{pade_approx}),(\ref{main_eq}) for $n=3$ are chosen the following: ~$\gamma_1=0.972926132$, $\gamma_2=0.744418059$, $\gamma_3=0.150843924$, $\beta_1=0.004210420$, $\beta_2=0.081312882$, $\beta_3=0.414236605$, for which, as shown in \cite{Lee01101985,seis_interpr}, such approximation is valid up to the angles of $89$ degrees. It is also possible to select large values $n$, thus asymptotically approaching to $\pi/2$, however in this case computer costs will essentially increase.

The method in question can be considered as analog of the spectral-difference one based on Fourier transform, but the part of "frequency"\; is played by the parameter $m$, which determines the degree of the Laguerre polynomials. An important property is in that opposed to the Fourier method the operator of problem (\ref{main_eq_laguerre}) is independent of the parameter separating harmonics $m$. Hence it appears possible to solve in the domain $(x,z,m)$ the common system of equations for different right-hand sides whereas this property is not valid in the domain $(x,z,\omega)$. After carrying out the spatial approximation of equations (\ref{main_eq_laguerre}), the multiple solution to the SLAE with a common matrix allows one to construct an effective computing procedure for solving a difference problem.

A disadvantage of the Laguerre transform is the absence of the fast transformation algorithm. However, taking into account the fact that input data are set only along the upper surface ($z=0$), and the inverse transformation is done for a fixed time instant, the total cost of the direct and inverse transformations  appears to be minor as compared to that needed for calculation of coefficients of decomposition (\ref{series_lag}) from the solution to problem (\ref{main_eq_laguerre}).
\subsection{Spatial approximation}
To approximate equation (\ref{main_eq1_laguerre}) we will use unconditionally stable the Crank-Nicolson scheme\cite{CrankNicolson2} of second order of accuracy:
\begin{subequations}
\label{main_eq_diff}
\begin{empheq}[left=\empheqlbrace]{align}
c\frac{ \bar{u}^m_{ik+1}-\bar{u}^m_{ik}}{h_z}+\tilde{\eta}\frac{\bar{u}^m_{ik}+\bar{u}^m_{ik+1}}{2}-\tilde{\eta}\sum_{s=1}^{3}\bar{\psi}^{m,s}_{ik+1/2}=
-\Phi_1\left(\frac{\bar{u}^m_{ik}+\bar{u}^m_{ik+1}}{2}\right)+\sum_{s=1}^{3}\Phi_1\left(\bar{\psi}^{m,s}_{ik+1/2}\right)
,\\
c^2\gamma_s \mathcal{L}_{x}\bar{\psi}^{m,s}_{ik+1/2}-\tilde{ \eta}^2\bar{\psi}^{m,s}_{ik+1/2}+1/2c^2\beta_s \mathcal{L}_{x}\bar{u}^m_{ik+1}=-1/2c^2\beta_s \mathcal{L}_{x}\bar{u}^m_{ik}+\Phi_2(\bar{\psi}^{m,s}_{ik+1/2}),\quad s=1,2,3,
\end{empheq}
\end{subequations}
where the difference operator $\mathcal{L}_{x}$ is of the form
\begin{equation}
\mathcal{L}_{x}f(x)\equiv\frac{1}{h_x^2}\left[a_0f(x)+\sum_{j=1}^{N}a_j\left((f(x-jh_x)+f(x+jh_x)\right)\right]=\frac{\partial^2 f}{\partial x^2}(x)+O(h_x^{ 2N}).
\label{drp_coeff}
\end{equation}
Formally, the presence of a singular component in the solution \cite{Bamberger1988,Bamberger1988a} does not allow attaining a high order of convergence for finite difference schemes, but for practical purposes the calculations have shown (see Section Numerical Experiments) that schemes of high orders of accuracy make possible to obtain solutions with a considerably lesser number of mesh nodes because such schemes more precisely reproduce the dispersion law. In deriving equations (\ref{main_eq}) the velocity model was assumed to be homogeneous, although it also yields satisfactory results for inhomogeneous media. In the latter case this model correctly keeps kinematics of waves, but not their amplitudes. For a wide range of problems such an approximate one-way model is admissible, as the correct account of amplitudes essentially increases computer costs \cite{Zhang2005}. Based on the above-said for approximating $\partial^2/\partial x^2$  it is reasonable to use the dispersion-relationship-preserving method (DRP) by Tam and Webb \cite{Tam1993262}, whose main idea is in the following. According to the Fourier derivative rule, $k_j\Longleftrightarrow-\mathrm{i}\partial_j,$
values of optimized coefficients $a_n$ in (\ref{drp_coeff}) are defined as solution to the problem of minimizing the error functional in the space of wave numbers
$$
\frac{\partial}{\partial a_n}\int_0^{k_c}\left|-k_x^2+\frac{1}{h_x^2}\left[a_0+\sum_{j=1}^{N}a_n\cos(jk_xh_x)\right]\right|^2w(k_x)dk_x=0,
$$
where $w(k_x)$ is some weight function, $k_c$ is a maximum accurate wave number. This approach and its various modifications \cite{GPR:GPR12160,Chu2012,Zhang2013511} make possible to decrease the number of mesh nodes and to preserve high accuracy of calculations as compared to conventional difference schemes obtained with the Taylor expansion in series \cite{Samarskii2001}.

\subsection{Solution of the SLAEs}
Solving the difference equations as a stage of the algorithm proposed demands high computer costs, therefore its applicability depends on the speed of solving the corresponding SLAEs. As opposed to Fourier transform, the coefficients of the Laguerre expansion in series (\ref{series_lag}) are dependent in a recurrent manner (see (\ref{main_eq_laguerre})). Hence, for a fixed $k$ for different $m$ it is required to solve SLAEs many times with a common real matrix and different right-hand sides. This allows the use of the methods based on $LU$ decomposition, where factorization of a matrix for all $m$ is done only once.

Let us write down the difference problem  (\ref{main_eq_diff}) in the form of a SLAE as
\begin{equation}
\label{main_slae}
\begin{array}{l}
\left(
\begin{array}{llll}
\gamma_1 c^2\mathcal{L}_{x}-\tilde{\eta}^2I&0&0&1/2\beta_1 c^2\mathcal{L}_{x}\\
0&\gamma_2 c^2 \mathcal{L}_{x}-\tilde{\eta}^2I&0&1/2\beta_2 c^2\mathcal{L}_{x}\\
0&0&\gamma_3 c^2\mathcal{L}_{x}-\tilde{\eta}^2I&1/2\beta_3 c^2\mathcal{L}_{x}\\
-2\tilde{\eta} I&-2\tilde{\eta} I&-2\tilde{\eta} I&\left(2c/h_z+\tilde{\eta}\right)I
\end{array}
\right)\left(
\begin{array}{l}
\bar{\mathbf{\Psi}}_{k+1/2}^{m,1}\\
\bar{\mathbf{\Psi}}_{k+1/2}^{m,2}\\
\bar{\mathbf{\Psi}}_{k+1/2}^{m,3}\\
\bar{\mathbf{U}}_{k+1}^m
\end{array}
\right) \\ = \left(
\begin{array}{c}
-c^2\beta_1\mathcal{L}_{x}\bar{\mathbf{U}}_{k}^m/2+\Phi_2\left(\bar{\mathbf{\Psi}}_{k+1/2}^{m,1}\right)\\
-c^2\beta_2\mathcal{L}_{x}\bar{\mathbf{U}}_{k}^m/2+\Phi_2\left(\bar{\mathbf{\Psi}}_{k+1/2}^{m,2}\right)\\
-c^2\beta_3\mathcal{L}_{x}\bar{\mathbf{U}}_{k}^m/2+\Phi_2\left(\bar{\mathbf{\Psi}}_{k+1/2}^{m,3}\right)\\
\left(2c/hz-\tilde{\eta}\right)\bar{\mathbf{U}}_{k}^m-\Phi_1\left(\bar{\mathbf{U}}_{k+1}^m+\bar{\mathbf{U}}_{k}^m\right)+2\sum_{s=1}^{3}\Phi_1\left(\bar{\mathbf{\Psi}}_{k+1/2}^{m,s}\right)
\end{array}
\right),
\end{array}
\end{equation}
where $I$ is the unit matrix. Employing the Schur complement \cite{Cottle1974}, the mesh functions $\bar{\mathbf{U}}_{k+1}^m$ can be calculated through the solution to the following reduced SLAE
\begin{equation}
\label{main_schur}
\begin{array}{l}
\left[\left(2c/h_z+\tilde{\eta}\right)I+\tilde{\eta}\sum_{s=1}^3\beta_s c^2 \mathcal{L}_{x}\left( \gamma_s c^2 \mathcal{L}_{x}-\tilde{\eta}^2I\right)^{-1}\right]\mathbf{\bar{U}}_{k+1}^m \\\\=\left(2c/hz-\tilde{\eta}\right)\bar{\mathbf{U}}_{k}^m-\Phi_1\left(\bar{\mathbf{U}}_{k+1}^m+\bar{\mathbf{U}}_{k}^m\right)+2\sum_{s=1}^{3}\Phi_1\left(\bar{\mathbf{\Psi}}_{k+1/2}^{m,s}\right)\\\\
+2\tilde{\eta}\sum_{s=1}^3\left( \gamma_sc^2 \mathcal{L}_{x}-\tilde{\eta}^2I\right)^{-1}\left(-c^2\beta_s\mathcal{L}_{x}\bar{\mathbf{U}}_{k}^m/2+\Phi_2\left(\bar{\mathbf{\Psi}}_{k+1/2}^{m,s}\right)\right),
\end{array}
\end{equation}
where the functions $\bar{\mathbf{\Psi}}^{m,s}_{k+1/2}$ are defined as
\begin{equation}
M_s \bar{\mathbf{\Psi}}^{m,s}_{k+1/2}=\tilde{\eta}^{-2}\left(-\beta_sc^2\mathcal{L}_{x}\left(\bar{\mathbf{U}}^m_{k}+\bar{\mathbf{U}}^m_{k+1}\right)/2+\Phi_2\left(\bar{\mathbf{\Psi}}^{m,s}_{k+1/2}\right)\right),\; s=1,2,3.
\label{psi_non_reduce}
\end{equation}
Making use of the matrix property \cite{Henderson1981}
\begin{equation}
\label{matrix:prop1}
\left(B+I\right)^{-1}B=I-\left(B+I\right)^{-1},
\end{equation}
we can make a simplification
\begin{equation}
\label{schur2}
\begin{array}{l}
\left[
 \left(2c/h_z+\tilde{\eta}+\tilde{\eta}\sum_{s=1}^3\frac{\beta_s}{\gamma_s}\right)I+\tilde{\eta}\sum_{s=1}^3\frac{\beta_s}{\gamma_s} M_s^{-1}\right]\mathbf{\bar{U}}_{k+1}^m  =\mathbf{\bar{F}}_u^m+\eta\sum_{s=1}^3M_s^{-1}\bar{\mathbf{F}}^m_{\psi_s} ,
 \end{array}
\end{equation}
where
$$
\begin{array}{ll}
\displaystyle
M_s=\frac{\gamma_s c^2}{\tilde{\eta}^2} \mathcal{L}_{x}-I,\\ \displaystyle
\mathbf{\bar{F}}_u^m=\left(2c/hz-\tilde{\eta}-\tilde{\eta}\sum_{s=1}^3\frac{\beta_s}{\gamma_s}\right)\bar{\mathbf{U}}_{k}^m-\Phi_1\left(\bar{\mathbf{U}}_{k+1}^m+\bar{\mathbf{U}}_{k}^m\right)+2\sum_{s=1}^{3}\Phi_1\left(\bar{\mathbf{\Psi}}_{k+1/2}^{m,s}\right),\\ \displaystyle
\bar{\mathbf{F}}^m_{\psi_s}=-\frac{\beta_s}{\gamma_s}\bar{\mathbf{U}}_{k}^m+\frac{2}{\eta^2}\Phi_2\left(\bar{\mathbf{\Psi}}_{k+1/2}^{m,s}\right).
\end{array}
$$
Multiplying equation (\ref{schur2}) by the matrix $M_1M_2M_3$ and taking into consideration the commutative property of \mbox{$M_iM_j=M_jM_i$}, we obtain the governing equation for the calculation of mesh functions $\mathbf{\bar{U}}^m_{k+1}$
\begin{equation}
\label{reduce_eq}
\begin{array}{ll}
\displaystyle
\left[M_1M_2M_3\left(2 c/h_z+\tilde{\eta}+\tilde{\eta}\sum_{s=1}^3\frac{\beta_s}{\gamma_s}\right)I+\frac{\tilde{\eta}\beta_1}{\gamma_1}M_2M_3+\frac{\tilde{\eta}\beta_2}{\gamma_2}M_1M_3+\frac{\tilde{\eta}\beta_3}{\gamma_3}M_1M_2\right]\mathbf{\bar{U}}^m_{k+1}\\\\
\displaystyle =M_1M_2M_3\mathbf{\bar{F}}_u^m+\tilde{\eta}\left(M_2M_3\bar{\mathbf{F}}^m_{\psi_1}+M_1M_3\bar{\mathbf{F}}^m_{\psi_2}+M_1M_2\bar{\mathbf{F}}^m_{\psi_3}\right).
\end{array}
\end{equation}
Matrix (\ref{reduce_eq}) is banded and can be explicitly represented without calculation of the matrices $M_s^{-1}$ thus allowing us to apply efficient algorithms for solving SLAEs based on $LU$-decomposition. Hence, we have to calculate $K-1$ different $LU$ decompositions for (\ref{reduce_eq}), where $K$ is the number of mesh nodes in the direction $z$. Note, also, that opposed to (\ref{schur2}), for calculating the vector in the right-hand side (\ref{reduce_eq}) the inversion of the matrix $M_s$ is not required.

After the calculation of the mesh functions $\bar{\mathbf{\mathbf{U}}}^{m}$ before turning to calculating the functions $\bar{\mathbf{\mathbf{U}}}^{m+1}$, the functions $\bar{\mathbf{\Psi}}^{m,s}$ can be expressed as
\begin{equation}
\label{psi_func}
\bar{\mathbf{\Psi}}^{m,s}_{k}=M_s^{-1}\left(-\frac{\beta_s}{\gamma_s}\bar{\mathbf{U}}^m_{k}+\frac{1}{\tilde{\eta}^2}\Phi_2\left(\bar{\mathbf{\Psi}}^{m,s}_{k}\right)\right)
-\frac{\beta_s}{2\gamma_s}\bar{\mathbf{U}}^m_{k},\quad s=1,2,3.
\end{equation}
Equation (\ref{psi_func}) is obtained applying the property of (\ref{matrix:prop1}) to (\ref{psi_non_reduce}).
The function $\bar{\mathbf{\Psi}}^{m,s}$ should be calculated at the integer mesh nodes, because this simplifies the implementation of the Richardson extrapolation algorithm considered in the next Section.
\subsection{Increasing accuracy of the spatial approximation via the Richardson extrapolation }
According to the studies \cite{fatab2011,Terekhov:2013,Terekhov2015206} it is needed to have a few tens of mesh nodes per a minimum wavelength in order to provide a reasonable accuracy of calculating schemes of second order for real spatial-temporal scales. In order to reduce computer costs, we increase the accuracy of the algorithm up to $O(h^4_z)$ based on the Richardson extrapolation \cite{rohtua,Marchuk:99314}.

Let the mesh functions $\bar{\mathbf{U}}(\omega_1)$, $\bar{\mathbf{U}}(\omega_2)$, defined on the meshes $\omega_1,\omega_{2}$ with the mesh steps $h_x,h_z$ and $h_x,h_z/2$,  be solutions to problem (\ref{main_eq_diff}).  Then the linear combination
\begin{equation}
\bar{\mathbf{U}}=\frac{1}{3}\left(4\bar{\mathbf{U}}(\omega_2)-\bar{\mathbf{U}}(\omega_1)\right)
\label{richardson1}
\end{equation}
approximates the solution to problem (\ref{main_eq_laguerre}) on the mesh $\omega_1$ accurate to $O(h_z^4)$. The Richardson extrapolation occurs more rarely than difference schemes of high accuracy orders. However, preliminary calculations have shown that multi-step methods of the Adams-Moulton and the Adams-Bashfort types do not provide stability even for small steps $h_z$. The implicit Runge-Kutta methods of high accuracy orders are not efficient because they demand multiple calculation of the right-hand side of the equation to be solved. The instability of methods of high orders is due to the presence of a singular component in the solution, but according to computational experiments the Richardson extrapolation makes possible to stabilize such instability.

{\bf Fourth order downward-continuation algorithm with the global correction }
To calculate the mesh functions $\bar{\mathbf{U}}^m_k,\bar{\mathbf{\Psi}}^m_k,\; k=1,...,K$ accurate to $O(h_x^{\xi}+h_z^4)$, the following is necessary:
 \begin{enumerate}
   \item Based on the cubic splines interpolate values of the functions $\Phi_1(\bar{\mathbf{U}}^m_{k})$, preset on the mesh $\omega_1$, into nodes $\omega_2$.
   \item Based on the cubic splines interpolate values of the functions $\Phi_1(\bar{\mathbf{\Psi}}^m_{k})$,$\Phi_2(\bar{\mathbf{\Psi}}^m_{k})$, preset on the mesh $\omega_1$, into semi-integer nodes $\omega_1,\omega_2$.
   \item On the mesh $\omega_1$, applying equation (\ref{reduce_eq}), calculate the solution $\bar{\mathbf{U}}^{m}_{k}(\omega_1)$.
   \item On the mesh $\omega_2$, applying equation (\ref{reduce_eq}), calculate the solution $\bar{\mathbf{U}}^{m}_{k}(\omega_2)$.
   \item Based on the Richardson extrapolation, correct the mesh function $\bar{\mathbf{U}}$ with (\ref{richardson1}).
   \item After the calculation of the mesh functions $\bar{\mathbf{\mathbf{U}}}^{m}_{k+1}$ for all depths $k=2,...,K$, the functions $\bar{\mathbf{\Psi}}_{k}^{m,s}$ can be defined from the solution to equations (\ref{psi_func}).
\item Turn to the calculation of $(m+1)\text{th},\;(m+2)\text{th}$, etc. coefficients of the expansion of the Laguerre series.
 \end{enumerate}

Also, we can formulate the second algorithm where the extrapolated value of (\ref{richardson1}) is calculated immediately after calculating the function $\bar{\mathbf{U}}^{m}_{k+1}(\omega_1)$, $\bar{\mathbf{U}}^{m}_{k+1/2}(\omega_2)$, $\bar{\mathbf{U}}^{m}_{k+1}(\omega_2)$, then the extrapolated value $\bar{\mathbf{U}}^{m}_{k+1}$  is used at the next step  for calculating $\bar{\mathbf{U}}^{m}_{k+2}(\omega_1)$, $\bar{\mathbf{U}}^{m}_{k+3/2}(\omega_2)$, $ \bar{\mathbf{U}}^{m}_{k+2}(\omega_2)$. This algorithm will be called "the downward-continuation algorithm with the local correction". If for calculating the solutions $\bar{\mathbf{U}}(\omega_1)$ and $\bar{\mathbf{U}}(\omega_2)$ one uses a numerically stable algorithm (the Crank-Nicolson scheme in our case), the extrapolation procedure based on the global correction will also be numerically stable, whereas in the general case, the stability is not ensured with the use of the local correction \cite{Marchuk:99314}.
\subsection{Parallel consideration}
The possibility of the efficient use of modern supercomputer systems is one of obligatory demands imposed upon the development of new numerical algorithms. Therefore let us consider specific features of the parallel implementation of the approach proposed. As opposed to Fourier transform, applying the Laguerre transform  does not allow the calculation in parallel of the function $\mathbf{U}^{m}$ for different $m$, as $(m+1)$th coefficient of series expansion (\ref{series_lag}) recurrently depends on the $m$-th coefficient (\ref{main_eq_laguerre}). Hence, it is needed to parallelize the algorithm at the stage of solving problem (\ref{main_eq_diff}). In this case, the main difficulties are in solving SLAEs (\ref{reduce_eq}) and  (\ref{psi_func}), whereas the parallel operation of multiplication by the matrices $M_1,M_2,M_3$ is not consuming from the communication standpoint. Let us consider the 2D data decomposition (Fig.~\ref{fig:parallel}a), where computing nodes contain processors having the shared memory access allowing implementation of in-node inter-process communications considerably faster than the communications among the nodes. For solving SLAEs (\ref{reduce_eq}),(\ref{psi_func}) with banded matrices as a parallel algorithm it is reasonable to use the parallel dichotomy algorithm, which was developed for tridiagonal matrices \cite{terekhov:Dichotomy} and block-tridiagonal matrices \cite{Terekhov:2013}. With respect to the number of arithmetical operations, the dichotomy algorithm is comparable with other available algorithms; however the time needed for inter-process communications is considerably less in the dichotomy algorithm as compared to other algorithms. This is because the implementation of the dichotomy process on a supercomputer reduces to calculating the sum of series for distributed data. The commutative and associative properties of addition enable a considerable reduction in the total computation time with the use of inter-processor interaction optimization algorithms. First, the parallel dichotomy algorithm was developed for solving SLAEs with the same tridiagonal matrix but different right-hand sides. In \cite{Terekhov2}, the dichotomy algorithm was applied to solving SLAEs with the Toeplitz tridiagonal matrices. It was shown that the Toeplitz tridiagonal matrices were able to effectively solve SLAEs both with one and several right-hand sides. In \cite{fatab2011,Terekhov2015206}, the dichotomy algorithm was applied to implement a spectral-difference method to calculate acoustic and elastic wave fields. The authors could effectively use from $2$ up to $8192$ processors for one calculation and to obtain a highly accurate numerical solution of the dynamic problem of elasticity theory. Thus, the dichotomy algorithm is a powerful instrument for solving SLAEs with tridiagonal and block-tridiagonal (banded) matrices.

The parallel dichotomy algorithm makes possible to use thousands of processors, however its efficiency depends on the number of simultaneously solved equations and the size of a SLAE matrix. However, such dependence is considerably weaker as compared to other parallel algorithms for solving tridiagonal and block-tridiagonal SLAEs. Due to dependence of the solution $\bar{\mathbf{U}}^{m}_{k+1}$ on $\bar{\mathbf{U}}^{m}_{k}$ it is impossible to solve problems (\ref{reduce_eq}),(\ref{psi_func}) for all depths in parallel, which would allow reducing communication costs at the expense of striping the inter-processor exchanges. Therefore the parallelization for $x$-direction should be done within one computational node profiting from the shared memory for the fast data exchange among processors thus allowing a considerable reduction of communication costs.

To increase the number of processors participating in one calculation it seems reasonable to carry out parallelization in $z$-direction. In this case calculations can be effectively implemented by the conveyor principle, i.e. as soon as at the node number $q$ the functions $\mathbf{\bar{U}}^{m},\mathbf{\bar{\Psi}}^{m}$  are calculated, their values and those of their second derivatives  $\partial^2/\partial z^2$ at the lower boundary of a sub-domain are communicated to the computational node number $q+1$. At the node with number $q$, the calculation of functions with numbers $m+1$ starts, while at the node $q+1$--functions with numbers $m$. The value of the second derivative at the boundary of two sub-domains is required for setting boundary conditions when constructing the cubic interpolation spline within the Richardson extrapolation procedure. Taking into account the fact that the number of terms of series (\ref{series_lag}) essentially exceeds the number of computational nodes, the load of such computational conveyor will be optimal. The conveyor communication costs will be minor because only one inter-node exchange for calculation of the functions $\mathbf{\bar{U}}^{m},\mathbf{\bar{\Psi}}^{m}$ for one value $m$ is required.
\section{Numerical experiments}
Let us discuss a few tests that would allow the evaluation of the quality of the solution to be obtained as compared with other known algorithms.
Numerical procedures were implemented in Fortran-90 using the MPI library. The calculation was performed on the "Lomonosov" ~supercomputer of the Moscow State University. The supercomputer comprises Intel Xeon X5570 fore-core processors operating at $2.93$ GHz in the Infiniband QDR communication environment. Each computational node contains two processors and 12 GB of RAM.

The coefficients of difference scheme (\ref{drp_coeff}) were chosen as follows\cite{Zhang2013511}: $a_0=-3.12513824$, $a_1=1.84108651$, $a_2=-0.35706478$, $a_3= 0.10185626$, $a_4=-0.02924772$, $a_5=0.00696837$, $a_6= -0.00102952$, to provide the twelfth approximation order. The calculations have shown that as compared to classical difference schemes of high accuracy orders based on the Taylor expansion, the DRP approach makes possible to decrease the number of nodes in the difference scheme approximately by the factor of two with comparable accuracy. Taking into consideration the fact that computer costs for $LU$ decomposition of matrices (\ref{reduce_eq}),(\ref{psi_func}) are proportional to the third degree of the matrix band width and forward and backward substitution of $LU$ decomposition have costs proportional to the second degree, the effect of applying the DRP schemes is significant. Also, this allows one to decrease the demands for a minimum required random access memory for storing the matrices $L$ and $U$.

\subsection{Impulse Response}
In the first test we illustrate analyzing the accuracy by impulse responses. For the calculation we used a homogeneous medium model with the speed $250\;m/s$ and the size $1.6\;km\times1\;km$, the mesh steps being  $h_x=h_z=1\;m$ and $h_x=h_z=0.5\;m$. A point source was located at the center of the upper surface; the time dependence was as follows:
\begin{equation}
f(t)=\exp\left[-\frac{(2\pi f_0(t-t_0))^2}{g^2}\right]\sin(2\pi f_0(t-t_0)),
\label{impulse_source}
\end{equation}
where $f_0=30\mathrm{Hz},\;t_0=0.2s,\;g=4$. As compared to the Fourier transform, where the basis functions are uniquely defined, the two parameters, $\alpha$ and $\eta$, should be set for using the Laguerre transform (\ref{series_lag}). These parameters were experimentally chosen based on the analysis of the convergence rate of the Fourier--Laguerre series for the shifted function $f(t)$ with $t_0=T$, where $T$ is the upper boundary of the time interval for which the wave field is calculated. The parameters $\alpha$ and $\eta$ are chosen such that the function $f(t)$ with $t_0=T$ in the mean-quadratic norm is approximated accurate to $\varepsilon < 10^{-3}$. It should be regarded that as the value of $\eta$ grows, the values of the spatial derivatives of the functions $\bar{\mathbf{U}}^m,\bar{\mathbf{\Psi}}^m$ will also increase. For this reason, the size of the mesh has to be simultaneously decreased in order to prevent growth of the error of the spatial approximation as the parameter $\eta$ grows. The number of addends in series (\ref{series_lag}) was $n=2000$ for $T=3\;s$; the expansion parameters were $\alpha=0$ and $\eta=800$.

\begin{figure}[!h]
\centering

 \begin{subfigure}[b]{0.5\textwidth}
                        \includegraphics[width=\textwidth]{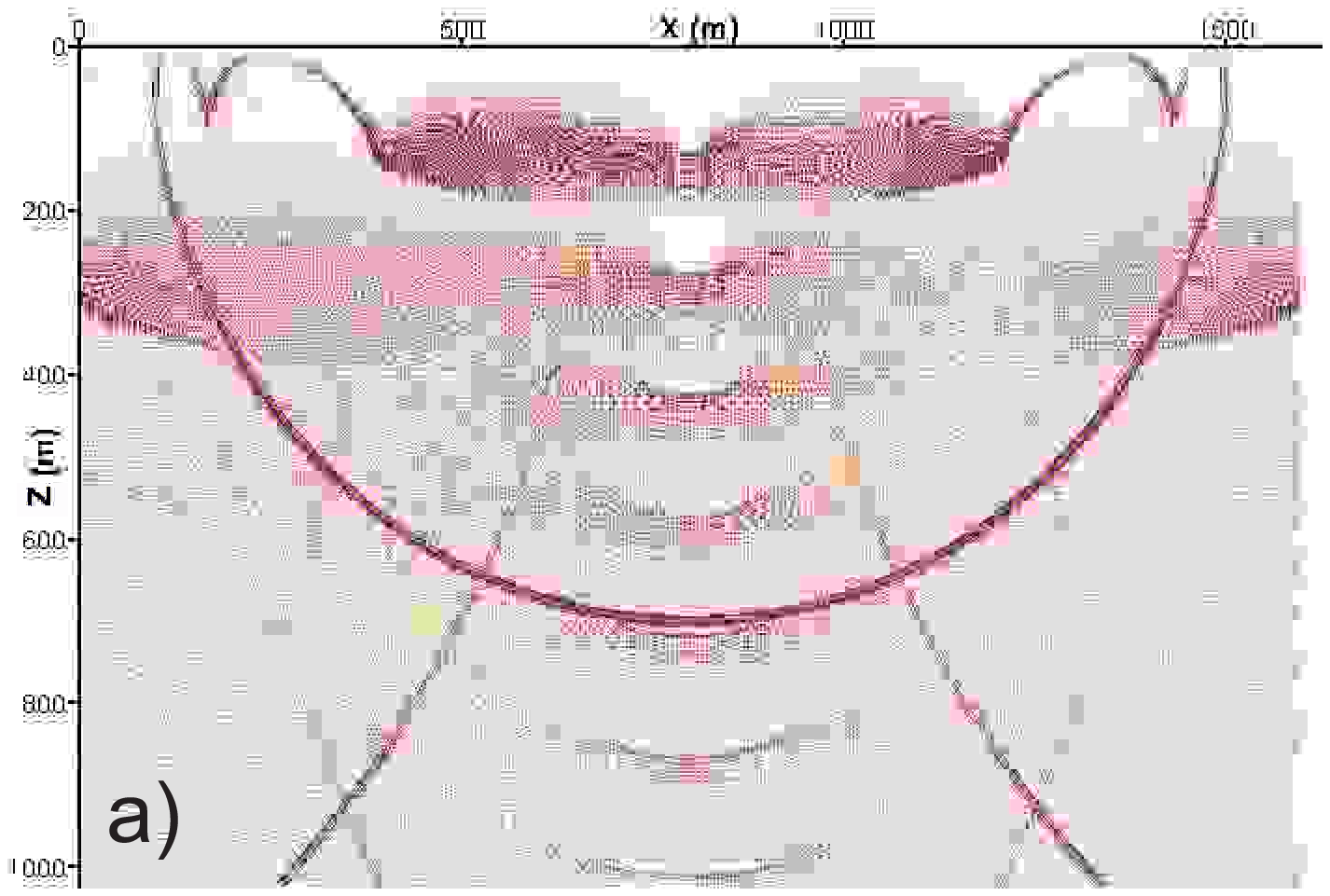}
                \label{fig:homogen1}
        \end{subfigure}%
        ~ 
                 \begin{subfigure}[b]{0.5\textwidth}
                \includegraphics[width=\textwidth]{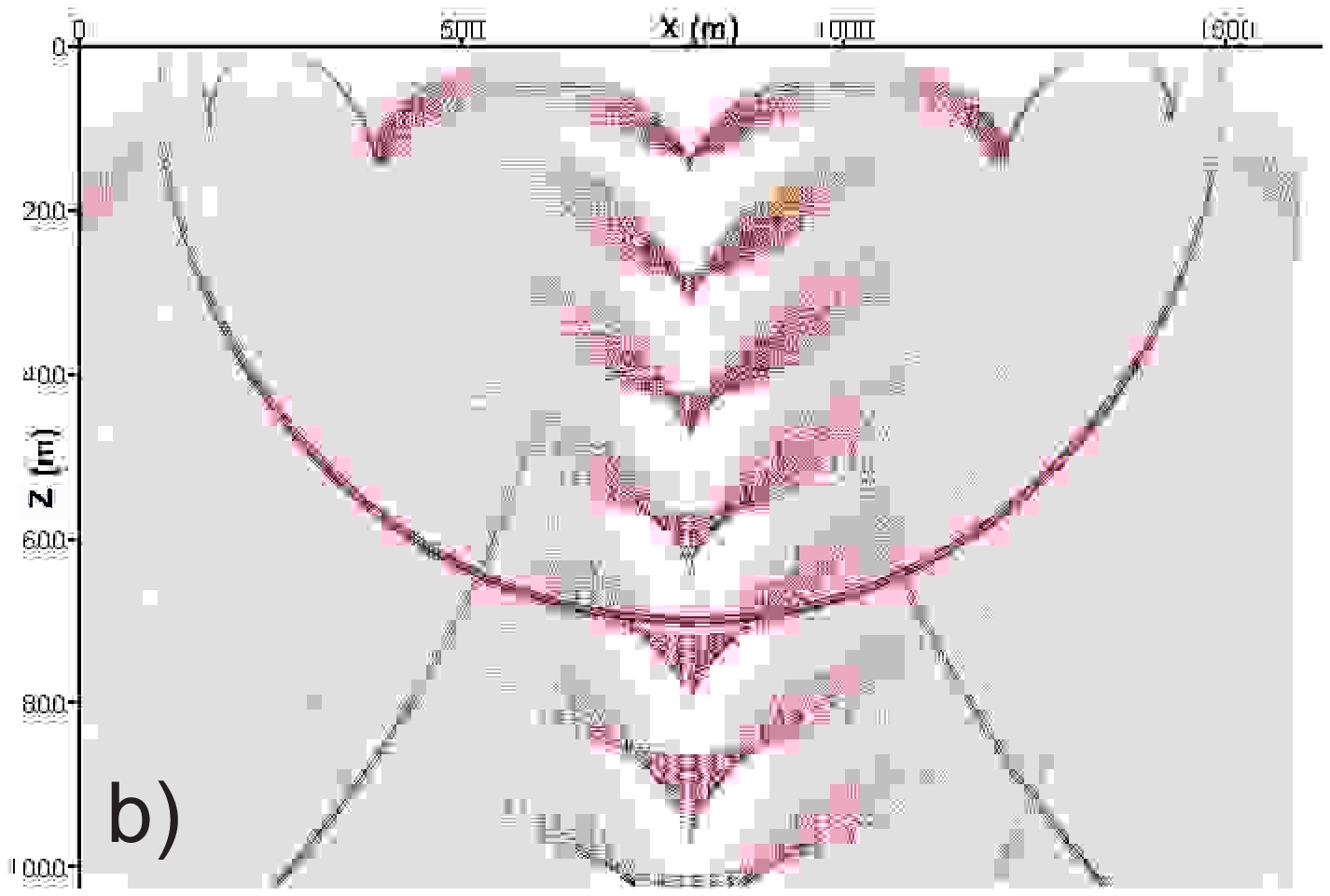}
                \label{fig:homogen2}
        \end{subfigure}%

        \begin{subfigure}[b]{0.49\textwidth}
                \includegraphics[width=\textwidth]{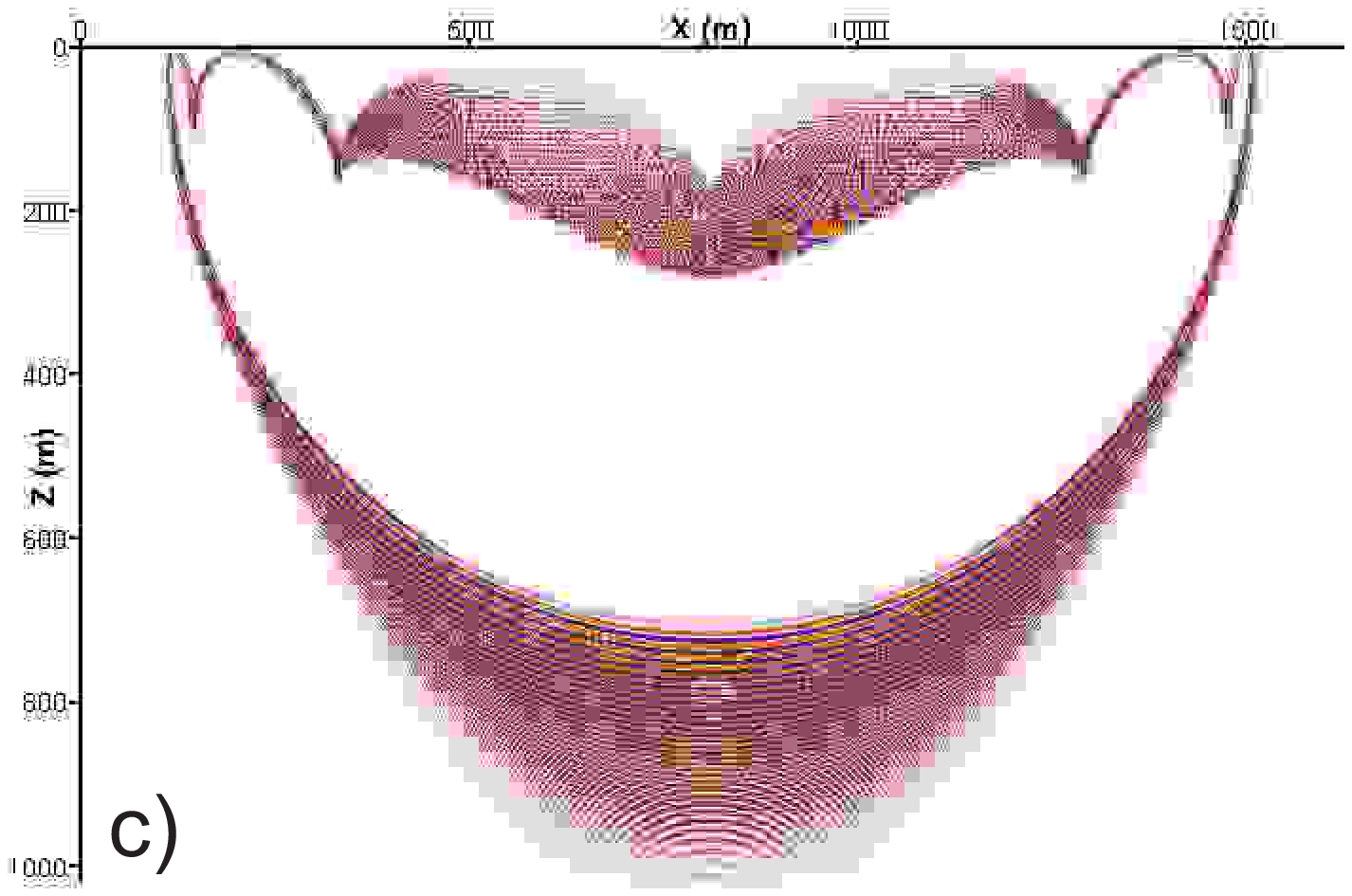}
                \label{fig:homogen3}
        \end{subfigure}
               ~ 
                \begin{subfigure}[b]{0.49\textwidth}
               \includegraphics[width=\textwidth]{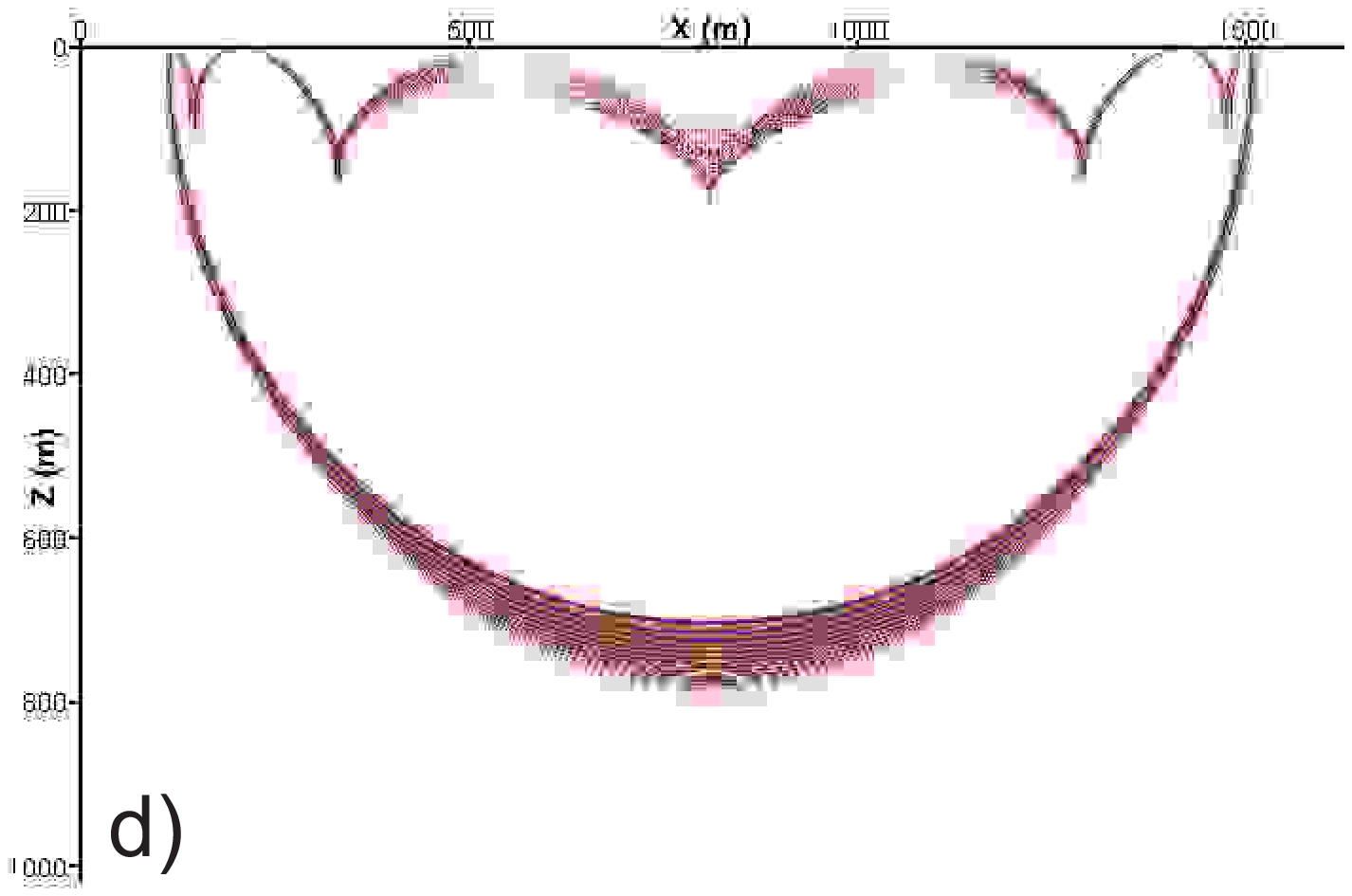}
                \label{fig:homogen4}
        \end{subfigure}%

        \begin{subfigure}[b]{0.49\textwidth}
                \includegraphics[width=\textwidth]{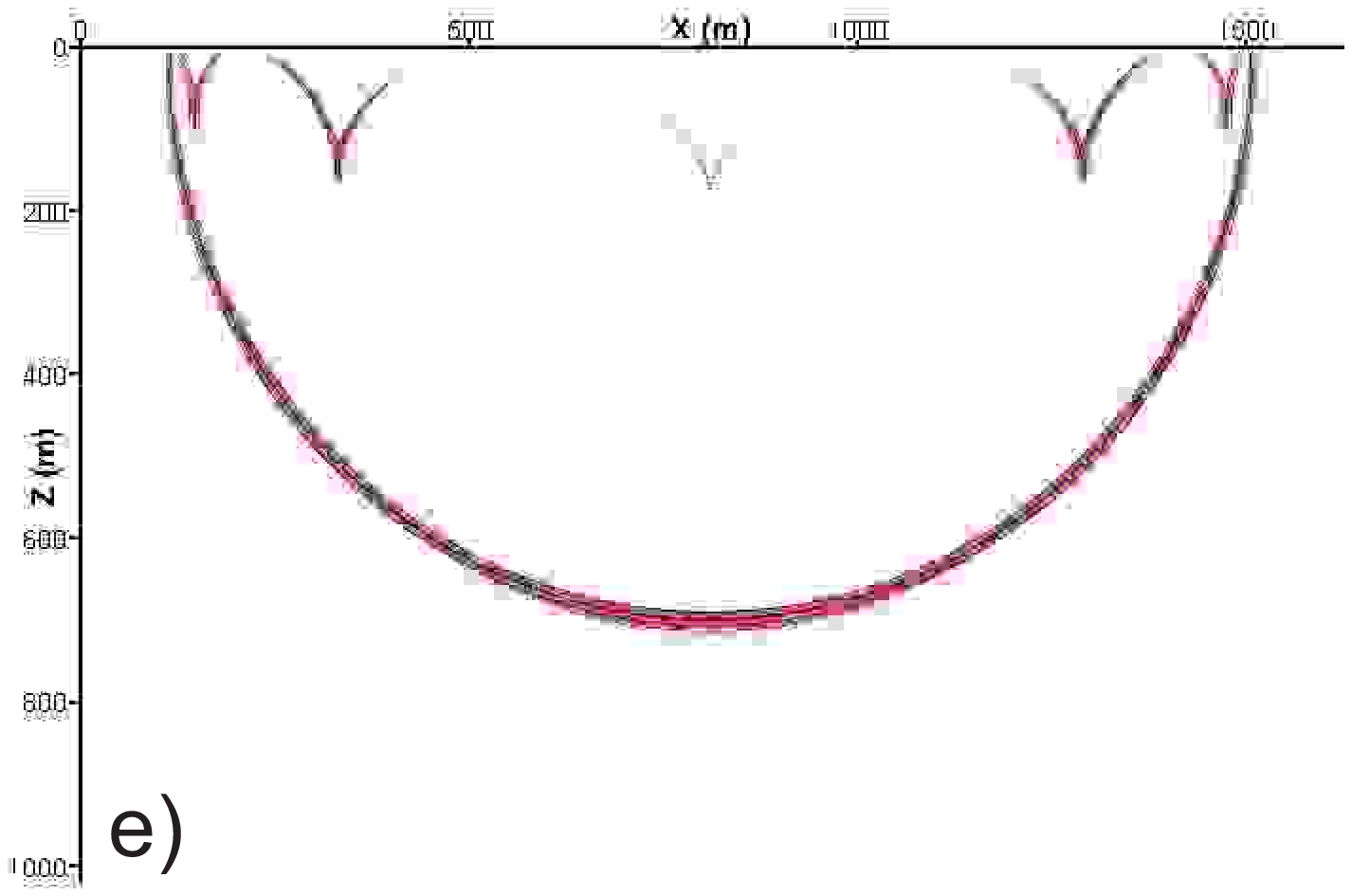}
                \label{fig:homogen5}
        \end{subfigure}
               ~
                \begin{subfigure}[b]{0.49\textwidth}
                \includegraphics[width=\textwidth]{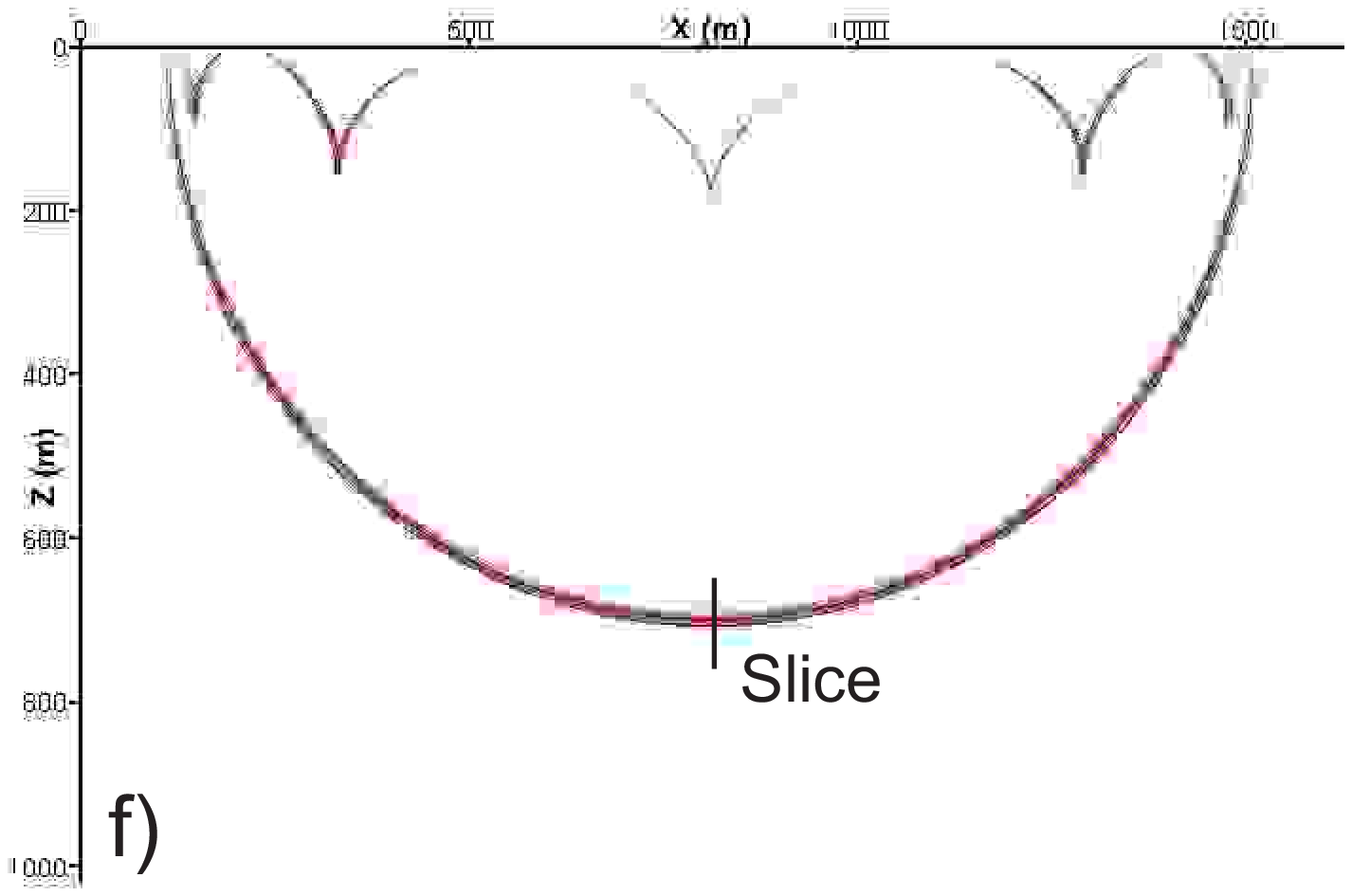}
                \label{fig:homogen6}
        \end{subfigure}%

\caption{Snapshots for the wave field at $t=3\;s$ for the homogeneous velocity model and impulse source $(\ref{impulse_source})$. FD method, 2nd order (a)~\mbox{$h_{x}=h_z=1\;m$} and (b)~\mbox{$h_{x}=h_z=1\;m$}, LFD method, 2nd order  (c)~\mbox{$h_{x}=h_z=1\;m$} and (d)~\mbox{$h_{x}=h_z=0.5\;m$}, LFD method, 4th order (e)~\mbox{$h_{x}=h_z=1\;m$} and (f)~\mbox{$h_{x}=h_z=0.5\;m.$}}
 \label{fig:homogen77}
\end{figure}

Figures~\ref{fig:homogen77}a,b present snapshots of the wave field from the point source for equation (\ref{pade_approx}), which were obtained employing the finite difference algorithm (FD) of second accuracy order using the Marchuk-Strang splitting. This software procedure has been implemented in the available program package Seismic Unix. It is evident (Figs.~\ref{fig:homogen77}a,b) that such an approach does not provide any accuracy of calculations due to multiple numerical noise and artifacts. In this case the transfer from the mesh with the step $1$ $m$ (Fig.~\ref{fig:homogen77}a) to the mesh with the step $0.5\;m$ (Fig.~\ref{fig:homogen77}b) does not allow an increase in calculation accuracy up to the acceptable level.

A conceptually other situation arises (Figs.~\ref{fig:homogen77}c,d,f,e) when the proposed Laguerre finite difference method (LFD) is used. It is clear that the new approach excludes the appearance of any numerical noise and the approximation error manifests as numerical dispersion which, as expected, is more distinct for the Crank-Nicolson scheme of second accuracy order (Figs.~\ref{fig:homogen77}c,d). In Figs.~\ref{fig:homogen77}e,f it is clear that applying the Richardson extrapolation makes possible to essentially improve the quality of solution. In addition, as for the second order scheme and for the fourth order of accuracy the transfer to a finer mesh allows diminishing the numerical dispersion. Thus, the key point in this case is the presence of convergence with decreasing the mesh step, which is not observed for calculations (Figs.~\ref{fig:homogen77}a,b).

\begin{figure}[!h]
\centering
        \begin{subfigure}[b]{0.49\textwidth}
                \includegraphics[width=\textwidth]{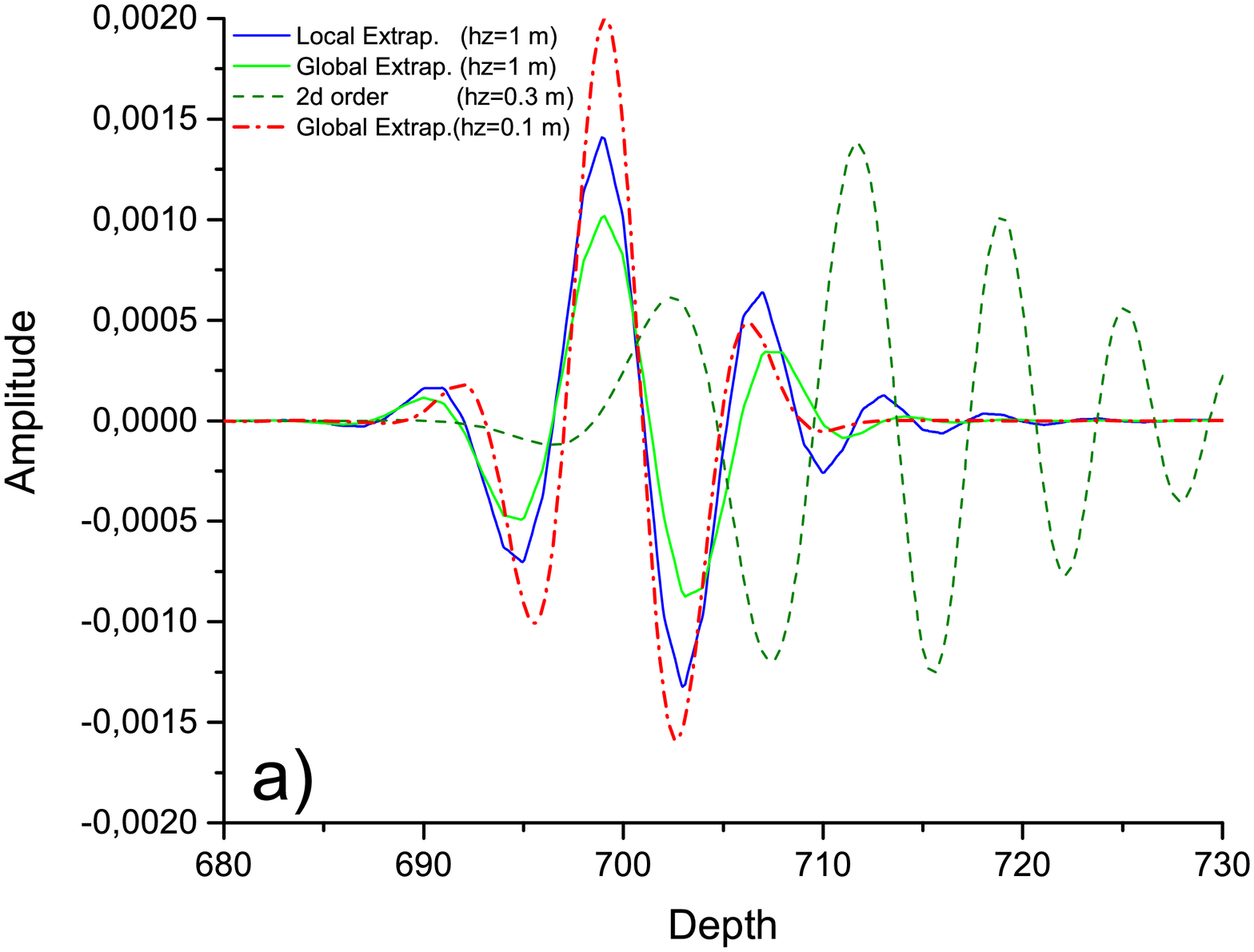}
        \end{subfigure}
               ~ 
                \begin{subfigure}[b]{0.49\textwidth}
                \includegraphics[width=\textwidth]{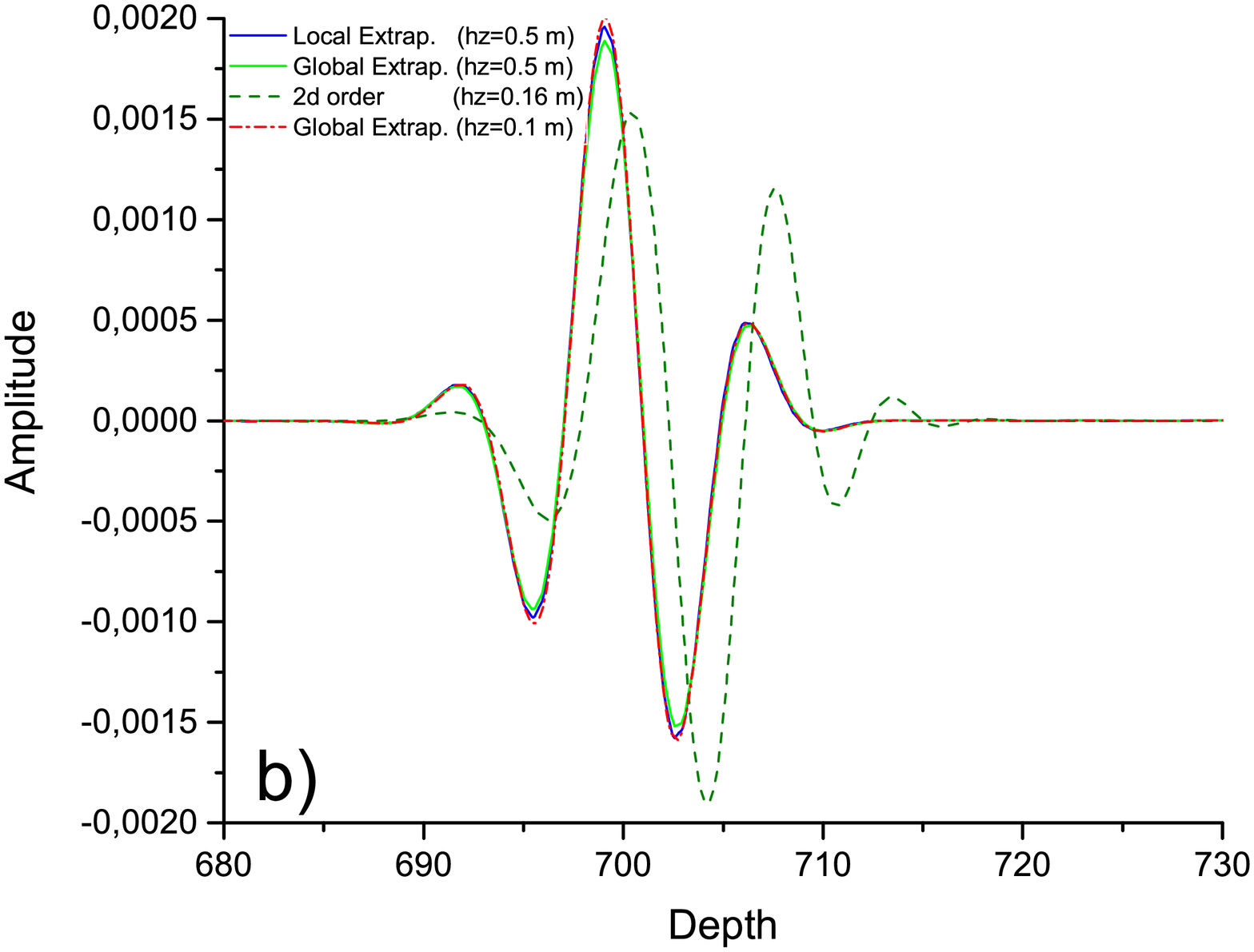}
        \end{subfigure}%

\caption{Dependence of the wave field on the coordinate along the straight line "Slice"\;(Fig.~\ref{fig:homogen77}) for different meshes for the LFD method.}
\label{fig:homo99}
\end{figure}

In spite of the absence of smoothness in the  solution due to singularities \cite{Bamberger1988,Bamberger1988a}, the use of the method of fourth order of accuracy allows a decrease in numerical dispersion as compared to the scheme of second order of accuracy. However one should take into account the fact that with an equal number of mesh points the algorithm based on the Richardson extrapolation requires three times as much arithmetical operations as the Crank-Nicolson scheme. This is explained by the necessity of solving equations (\ref{main_eq_diff}) on two meshes with the steps $h$ and $h/2$. The question arises whether it is possible just to increase the number of mesh nodes by the factor of three and to make use of the second order scheme to obtain a result compatible in accuracy with the fourth order scheme. Fig.~\ref{fig:homo99} shows the dependence of a wave field on the coordinate along the straight line "Slice"\;(Fig.~\ref{fig:homogen77}f). As is obvious, when the mesh step for the method of second accuracy order is three times smaller than that for the fourth order method (computer costs being approximately the same) the second order method yields less accurate results. Thus, it is reasonable to apply more complicated procedures based on the Richardson extrapolation, requiring large computer costs. It is also evident that the algorithm which is based on the local extrapolation inserts additional dispersion errors but better preserves amplitudes. Nevertheless, we are not going to use it because it is unstable for inhomogeneous media and inserts distinct artifacts associated with numerical dispersion into the image of a wave field.

When implementing the Richardson extrapolation algorithm we considered selecting the interpolation procedure and its influence on the solution to be solved. When using the linear interpolation the LFD algorithm demonstrated a very strong anisotropic dissipation in $z$-direction such that when passing a distance of several wavelengths the wave amplitude was tending to zero. The use of the parabolic or the cubic Lagrange interpolation somewhat improved the situation, but preservation of the amplitude was worse than for the cubic spline interpolation. Using the Lagrange interpolation of high orders as well as using splines of the orders exceeding the third caused instability of calculation. Thus, the choice of the cubic interpolation spline can be considered to be appropriate as other interpolation algorithms are either unstable or essentially dissipative.

Although the stability of the algorithm proposed has not been proved yet due to the difficulty in making the analytical analysis. In view of the following: asymmetry of matrices, the presence of singularities in a solution, applying the Richardson extrapolation procedure spline interpolation, etc. the algorithm in question in numerous experiments for different velocity models appeared  to be more stable provided $h_z\leq h_x$ and a local contrast of the velocity medium model is not quite distinct. A higher level of stability as compared to other algorithms is explained by the fact that approximation errors in the direction of the wave field extrapolation are dominantly of dissipative character.
\subsection{Migration procedure}
To form the image of the Earth's interior by the recorded seismic data one can use different algorithms of seismic migration \cite{seis_interpr,Claerbout:1985,Gray2001,Biondi2006} including those of the wave field downward continuation algorithms. Let us consider the post-stack migration based on the model of explosive boundaries \cite{Claerbout1971}. This model allows the evaluation of the algorithm proposed as compared to other known methods. As many studies have dealt with the problem in question, we will focus on its computational features. As there is no analytical solution, we will make a comparison by the FD method, with the numerical-analytical Fourier Finite Difference (FFD) method \cite{Ristow01121994} and the Phase-Shift Plus Interpolation method (PSPI) \cite{Gazdag1984} that are widespread in computational seismic prospecting. The two latter algorithms were developed instead of the FD method because it does not provide the required accuracy of calculation, and what is more important the computational stability when solving the one-way equation. However the FFD and PSPI methods approximately calculate the solution to the one-way equation for an inhomogeneous velocity model. The reason is that for the PSPI method the solution is obtained as a combination of several solutions for inhomogeneous models with different velocities.  In the context of the FFD method, the input operator is represented as sum of two operators, one of them corresponding to a homogeneous medium with some averaged in horizontal velocity (usually, a minimum velocity is used as the reference velocity). The other operator represents a varying velocity correction.  In the context of the Marchuk-Strang splitting procedure the first operator is inversed by analytical methods and the solution for correction is defined by finite difference methods of second order of accuracy. Although such a change is not similar to the original problem in terms of mathematics and is unstable for strongly contrast media, such an approach makes possible to considerably decrease computer costs. On the other hand, impossibility of evaluating the level of numerical errors on a sequence of imbedded meshes is an essential disadvantage of the FFD algorithm.
\subsubsection{Syncline model}
\begin{figure}[!h]
\centering
        \begin{subfigure}[b]{0.49\textwidth}
                \includegraphics[width=\textwidth]{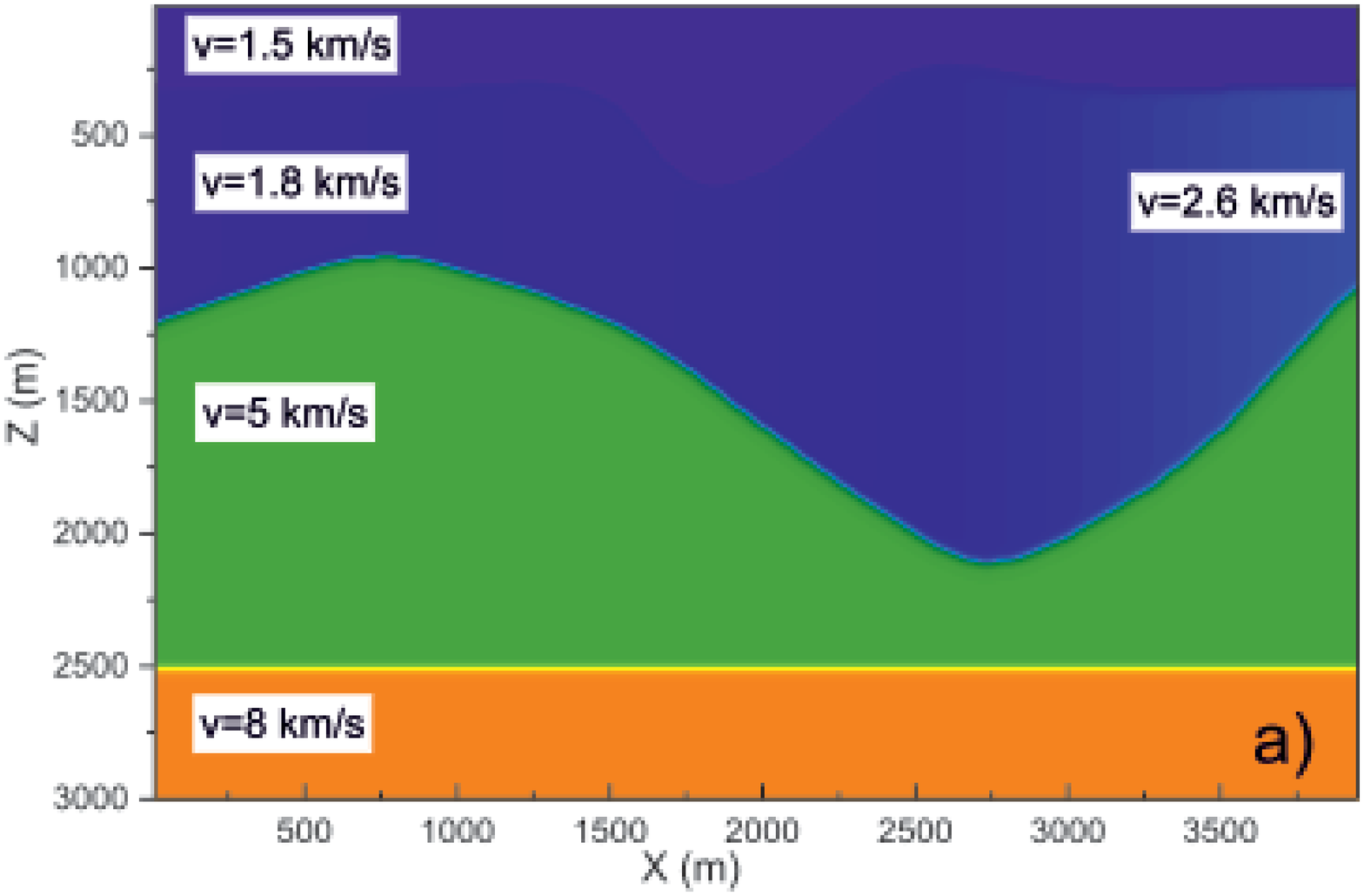}
        \end{subfigure}
              ~
                \begin{subfigure}[b]{0.49\textwidth}
                \includegraphics[width=\textwidth]{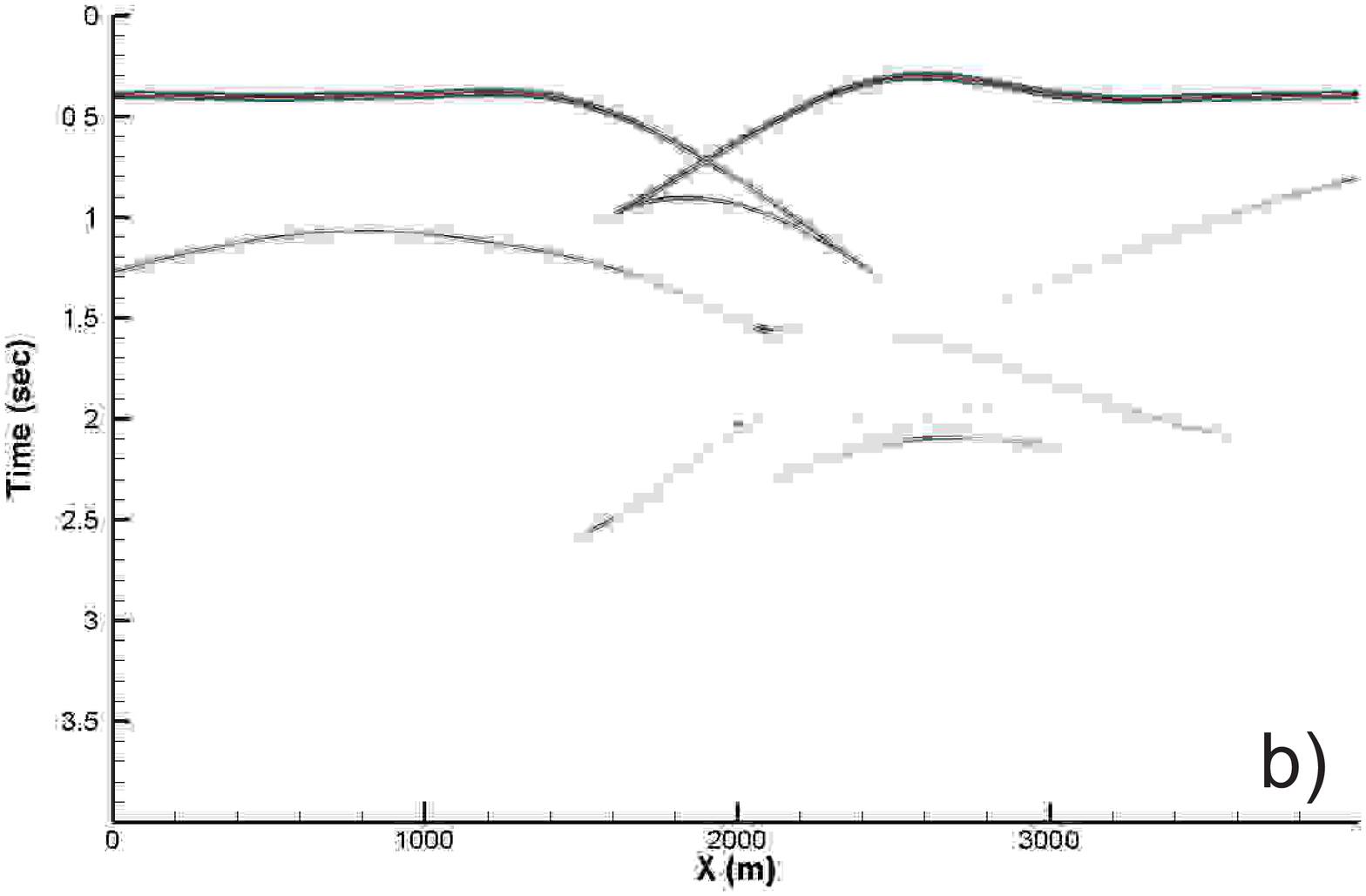}
        \end{subfigure}%
        \caption{(a) Syncline model and (b) zero-offset section.}
        \label{fig:syncline}
\end{figure}

Theoretical seismograms (Fig.~\ref{fig:syncline}b) for the syncline model (Fig.~\ref{fig:syncline}a) were obtained with the help of the Gaussian beams algorithm \cite{GaussianBeams_popov,Cerveny1985} implemented in the package Seismic Unix. For setting the boundary condition on the upper surface, the function for zero-offset section $\left.u(x,z,t)\right|_{z=0}=g(x,t)$. Was expanded in series (\ref{series_lag}) with the parameters $n=2000$, $\alpha=0$ and $\eta=600$ for $t\in[0,4]\;s$. The calculations were carried out on meshes with the steps $h_{x,z}=10\;m$, $5\;m$ and $0.5\;m$. According to the model of explosive boundaries the calculation velocities were set to be half the true velocity of the medium model.

\begin{figure}[!ht]
\centering
        \begin{subfigure}[b]{0.49\textwidth}
                \includegraphics[width=\textwidth]{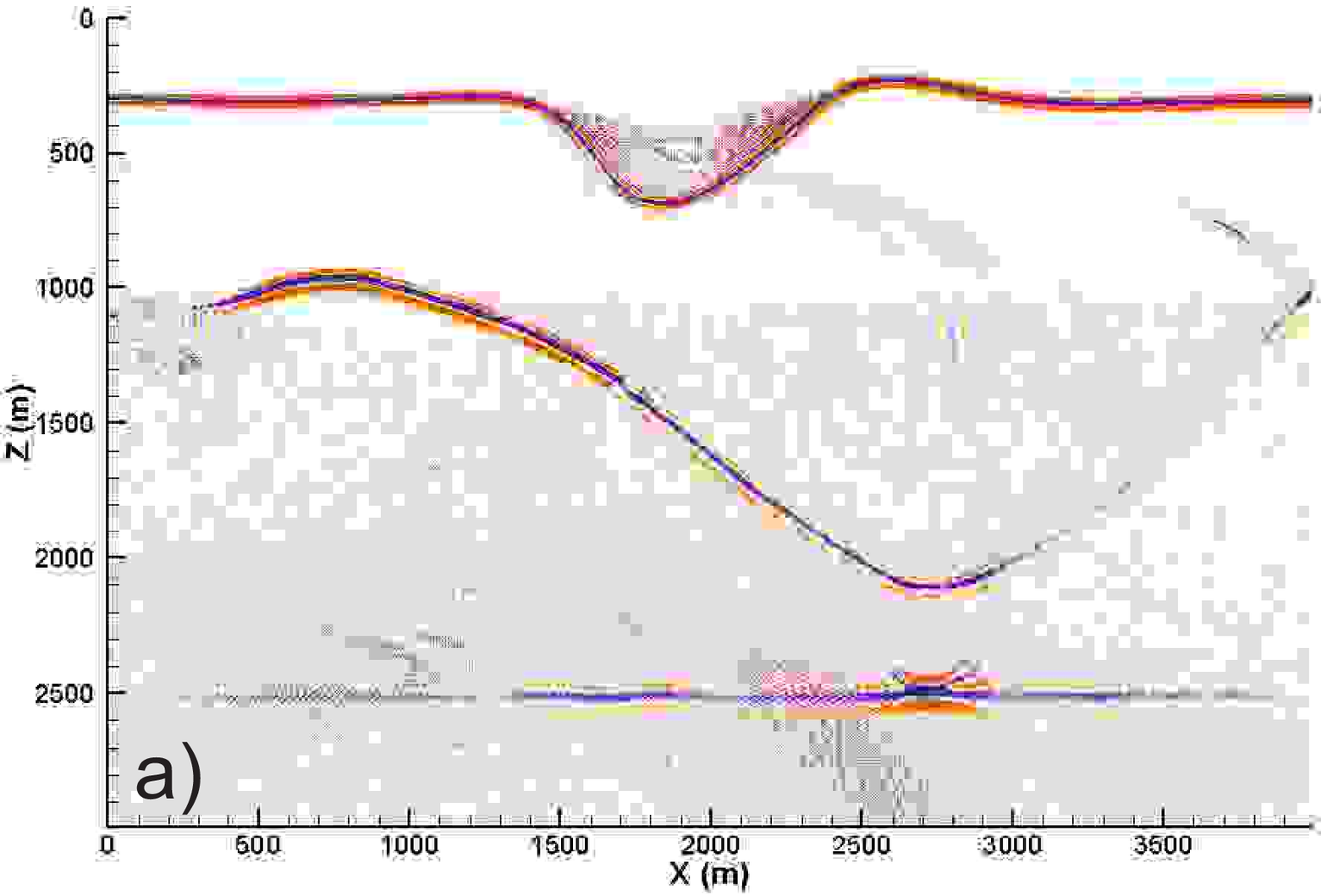}
                \label{fig:syn11}
        \end{subfigure}
               ~
        \begin{subfigure}[b]{0.49\textwidth}
                \includegraphics[width=\textwidth]{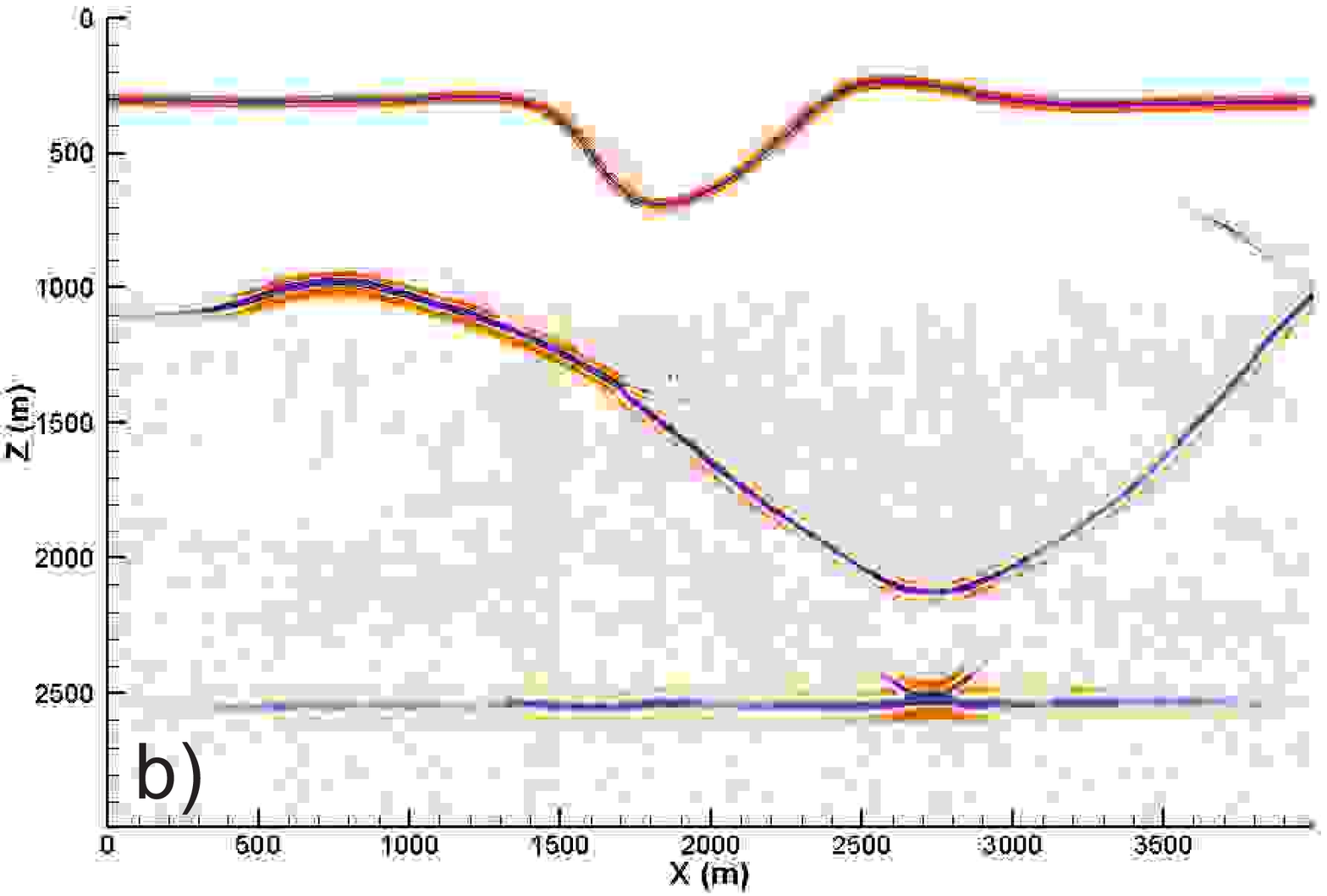}
                \label{fig:syn12}
        \end{subfigure}%
        \caption{Snapshots for the wave field at $t=4\;s$ for the model and the zero-offset section in Fig.~\ref{fig:syncline} for the FD method (a)~\mbox{$h_{x}=h_z=10\;m$} and (b)~\mbox{$h_{x}=h_z=5\;m.$} }
        \label{fig:syn1}
\end{figure}

\begin{figure}[!ht]
\centering
         \begin{subfigure}[b]{0.49\textwidth}
                \includegraphics[width=\textwidth]{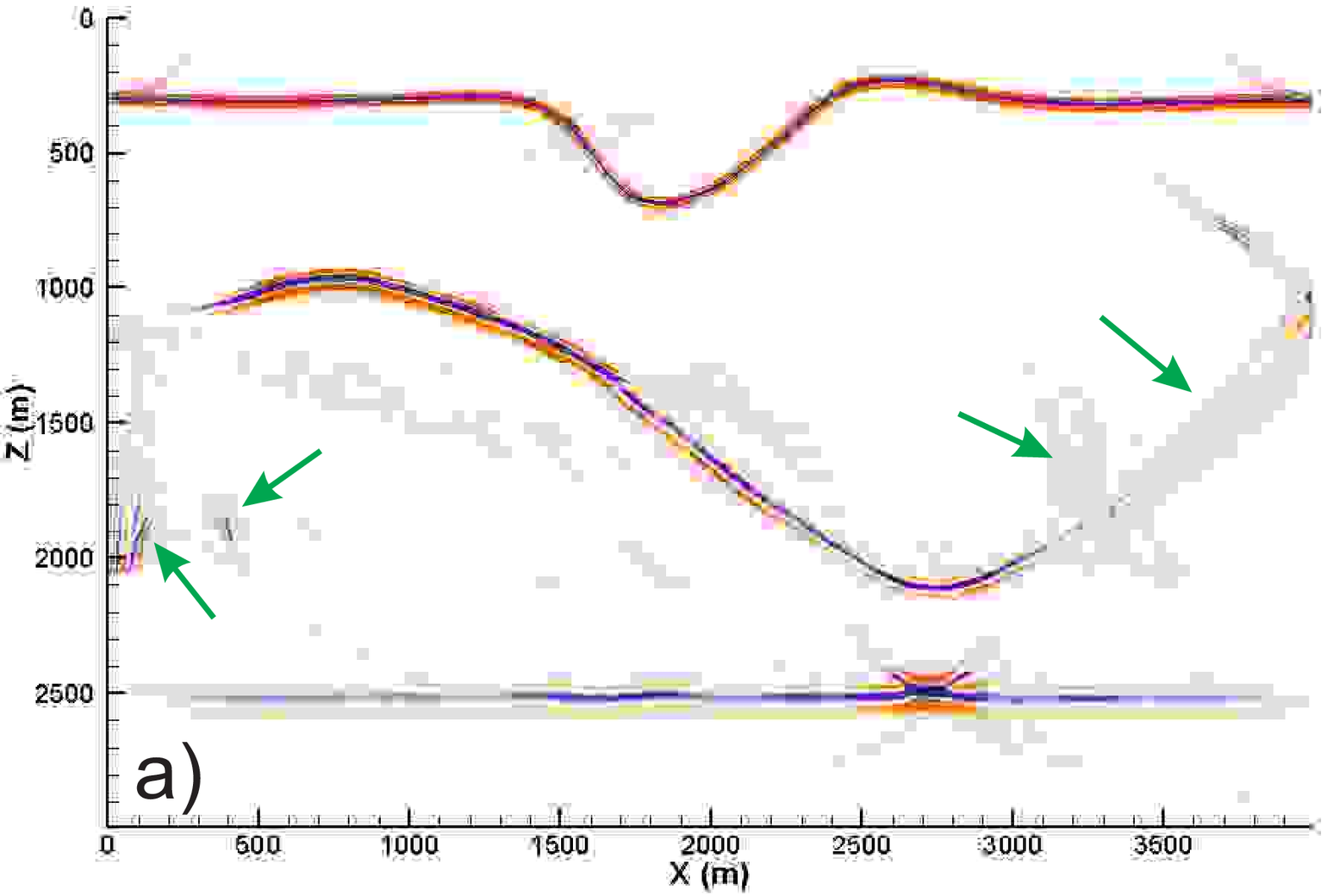}
                \label{fig:syn21}
        \end{subfigure}
               ~
                \begin{subfigure}[b]{0.49\textwidth}
                \includegraphics[width=\textwidth]{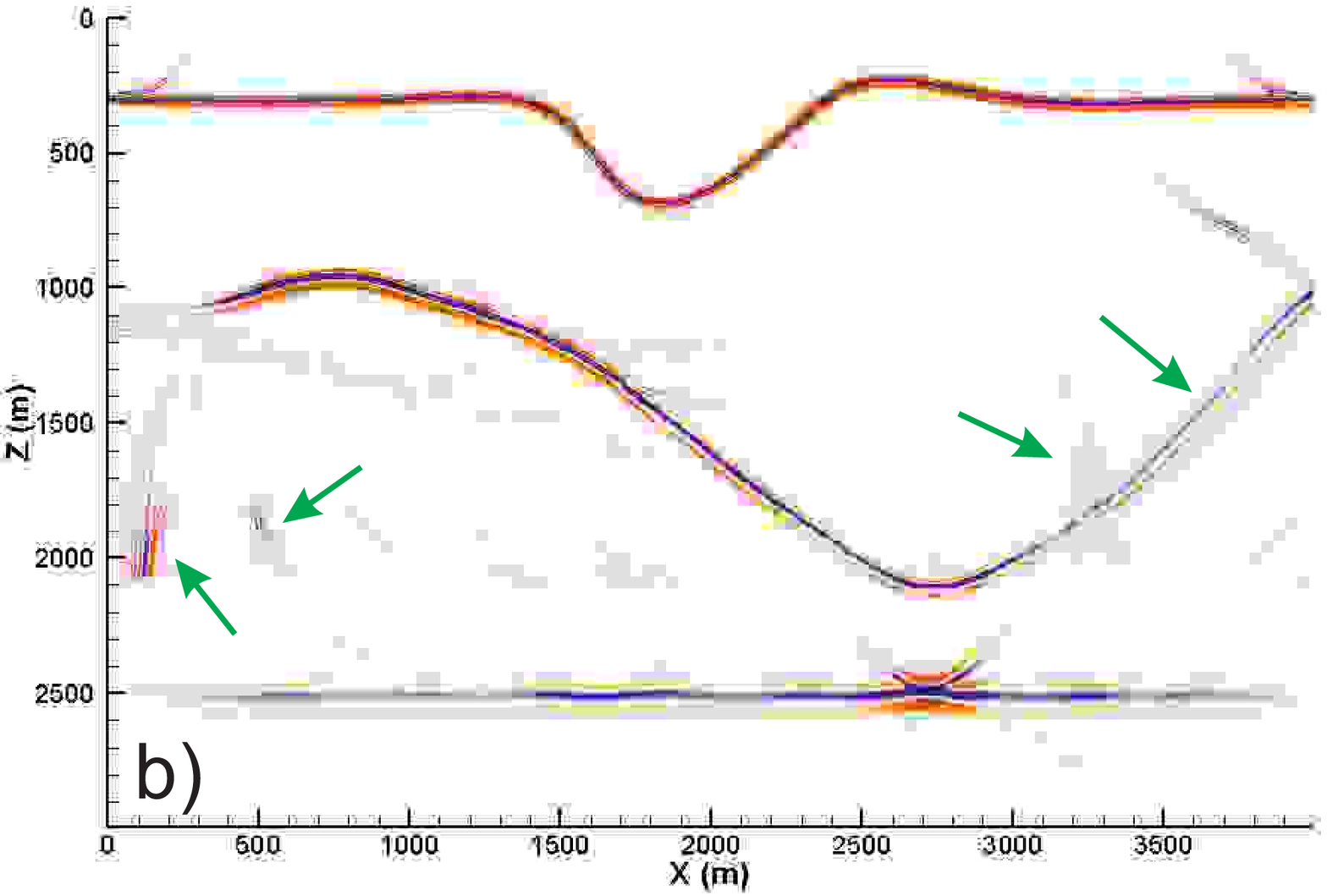}
                \label{fig:syn22}
        \end{subfigure}%
        \caption{Snapshots for the wave field at $t=4\;s$ for the model and the zero-offset section in Fig.~\ref{fig:syncline} for the FFD method (a)~\mbox{$h_{x}=h_z=10\;m$} and (b)~\mbox{$h_{x}=h_z=5\;m$}. Computational artifacts are indicated by arrows.}
        \label{fig:syn2}
\end{figure}

As well as in the calculation of the wave field from the point source (Fig.~\ref{fig:homogen77}) inadmissible level of numerical noise in the context of the migration problem for the FD method is observed (Fig.~\ref{fig:syn1}). In the FFD method, the noise level is much lower (Fig.~\ref{fig:syn2}), however, there are numerical artifacts and, in addition, the second horizon is insufficiently focused. When turning from the mesh with the step $h_z=10\;m$ (Fig.~\ref{fig:syn2}a) to the mesh with the step $h_z=5\;m$ (Fig.~\ref{fig:syn2}b) manifestation of artifacts remained as previously. Thus, the FFD algorithm can generate artifacts but does not allow their a posteriori evaluation on a sequence of imbedded meshes.

\begin{figure}[!h]
\centering
        \begin{subfigure}[b]{0.49\textwidth}
                \includegraphics[width=\textwidth]{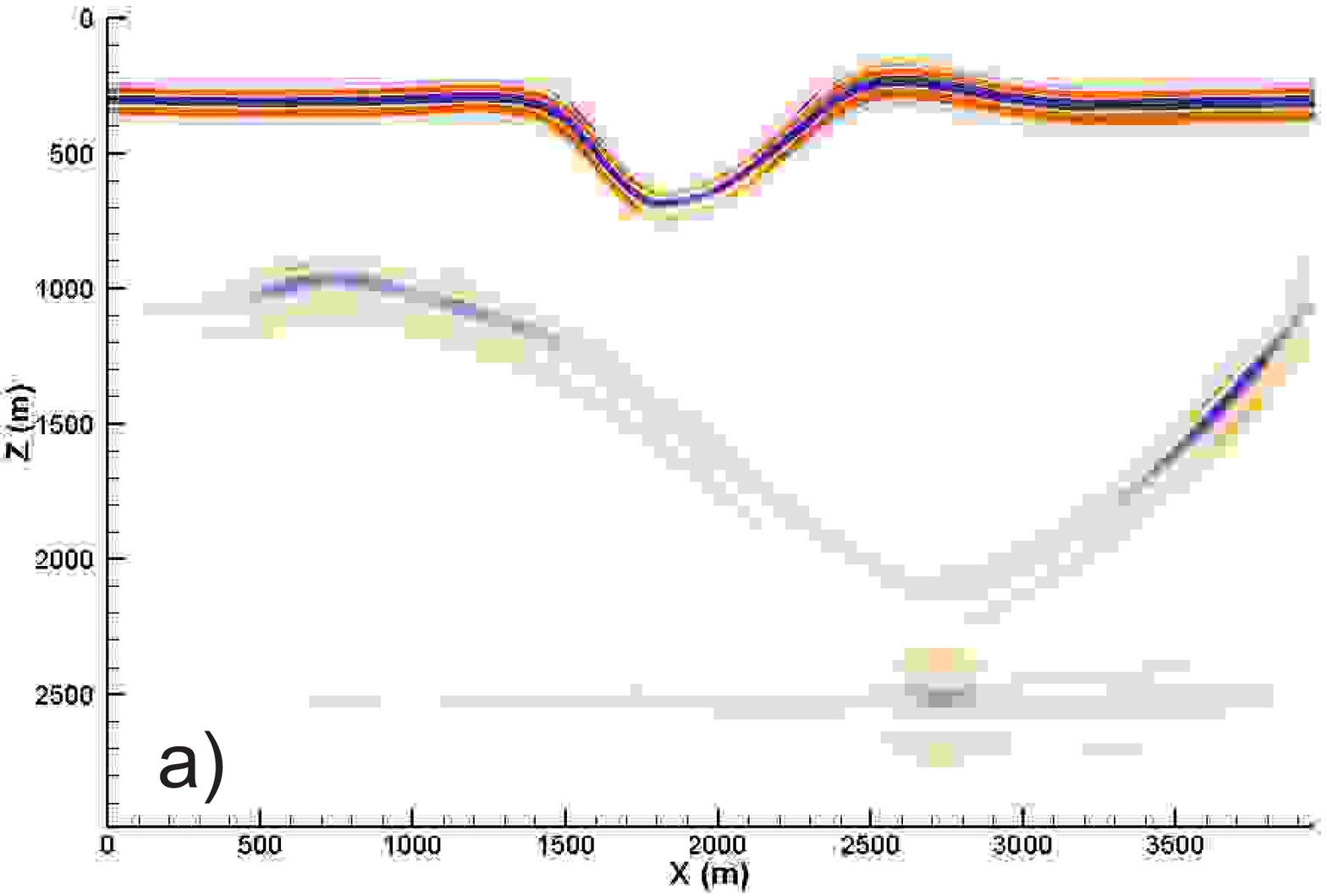}
        \end{subfigure}
               ~
                \begin{subfigure}[b]{0.49\textwidth}
                \includegraphics[width=\textwidth]{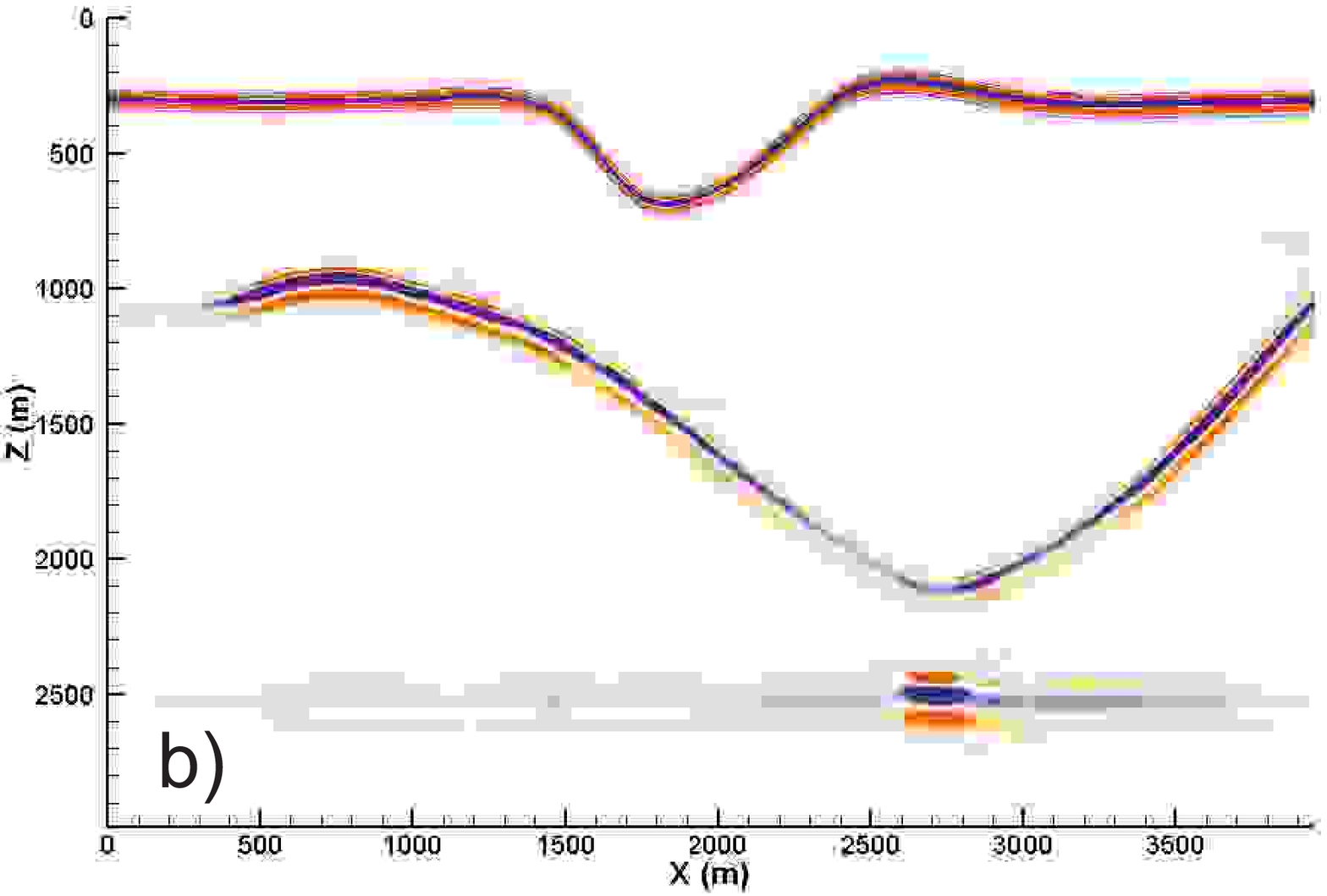}
        \end{subfigure}%

        \begin{subfigure}[b]{0.49\textwidth}
                \includegraphics[width=\textwidth]{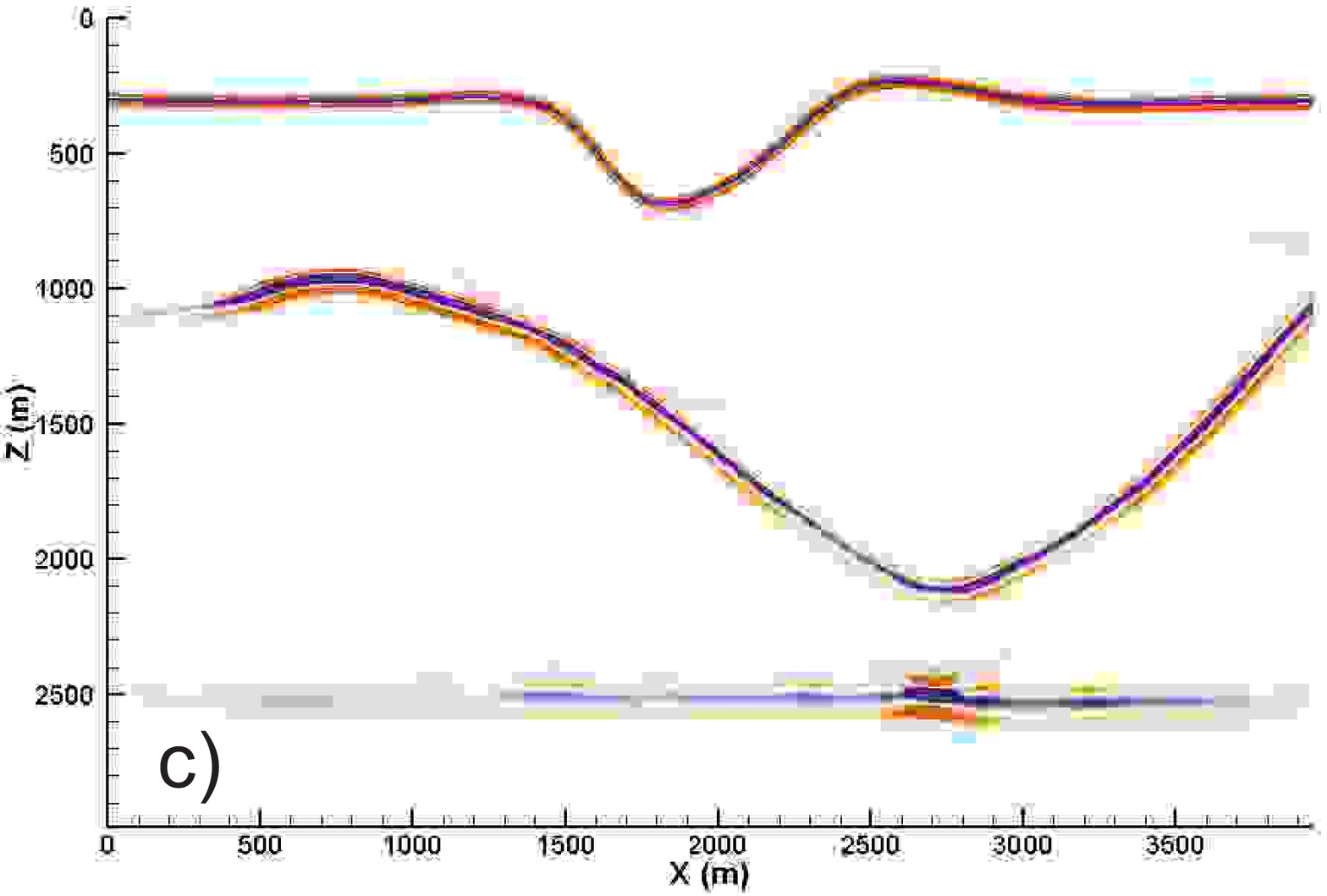}
        \end{subfigure}
        ~
        \begin{subfigure}[b]{0.49\textwidth}
                \includegraphics[width=\textwidth]{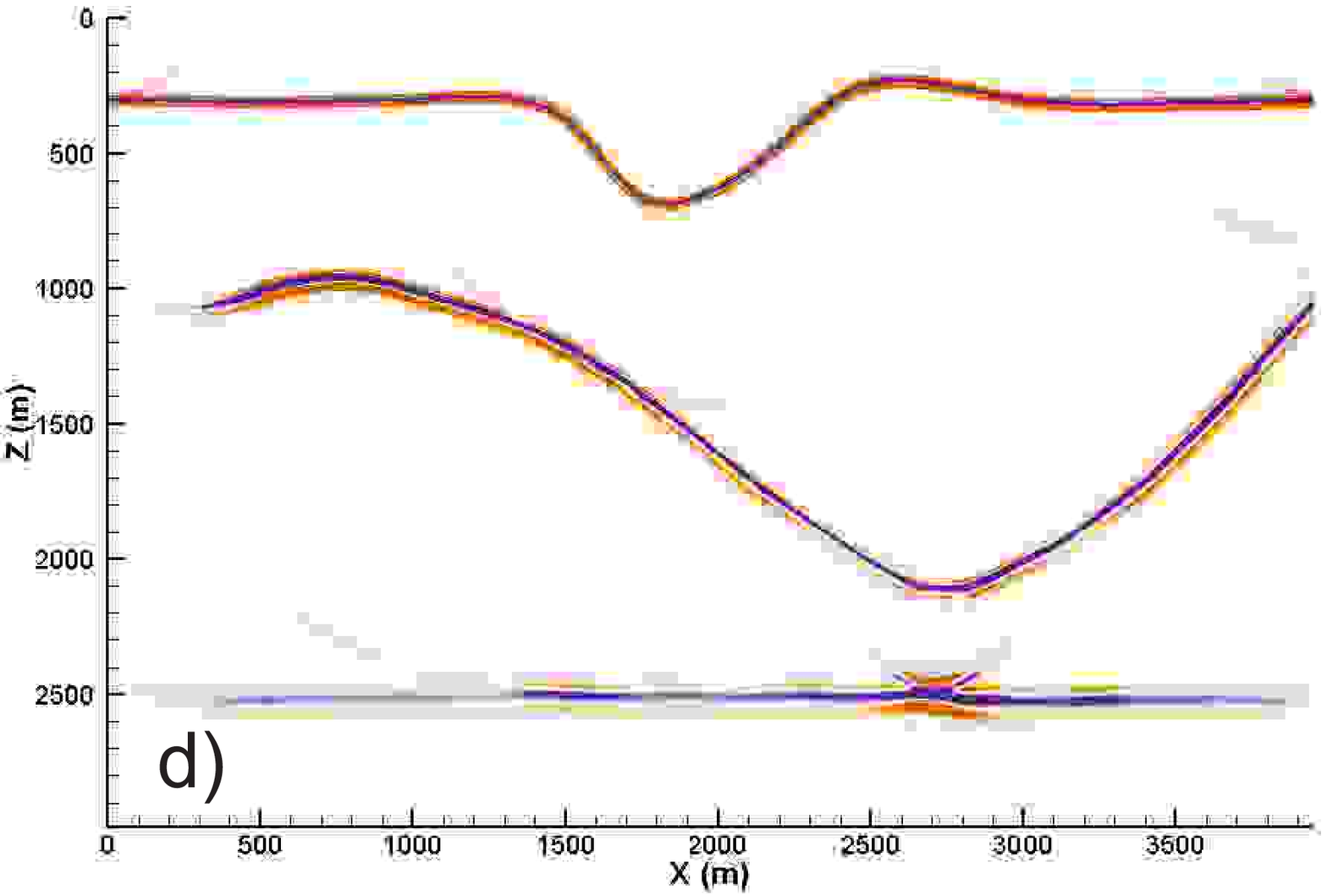}
        \end{subfigure}
        \caption{Snapshots for the wave field at $t=4\;s$ for the model and the zero-offset section in Fig.~\ref{fig:syncline} for the LFD method (a)~\mbox{$h_{x}=h_z=10\;m$}, (b)~\mbox{$h_{x}=10\;m,~h_z=5\;m$}, (c)~\mbox{$h_{x}=10\;m,~h_z=3\;m$}, (d)~$h_{x}=10\;m,$ $h_z=0.5\;m$}
        \label{fig:syn3}
\end{figure}

A conceptually different situation is observed with the LFD algorithm (Fig.~\ref{fig:syn3}), where the absence of high-frequency noise in the images obtained is distinct. The finite difference approximation errors for $x$-direction manifest as effect of numerical dispersion, whereas for $z$-direction the errors are of dissipative nature. The influence of these errors in turning to a finer mesh rapidly decreases, and the solution convergence of the mesh step makes possible to separate approximation errors from those of input data and inaccuracy of a model as it is. Hence, analyzing the results obtained when solving practical tasks can be essentially simplified. Additionally, calculations for a finer mesh with the step $h=0.5\;m$ for all the three algorithms were carried out. The FD method for a small mesh step appeared to be unstable, whereas the quality of the image  for the FFD method remained unchanged because the level of noise and artifacts was comparable with the calculation results with the mesh steps $h=5\;m$ and $h=3\;m$. For the LFD method, the image in Fig.~\ref{fig:syn3}d does not visually differ from that obtained with the mesh step $h=3\;m$ (Fig.~\ref{fig:syn3}c). Thus, with allowance for the convergence of the LFD algorithm, the result obtained with the mesh step $h=3\;m$ can be considered to be sufficiently accurate.

\subsubsection{Sigsbee Model}
\begin{figure}[!h]
\centering
        \begin{subfigure}[b]{0.49\textwidth}
                \includegraphics[width=\textwidth]{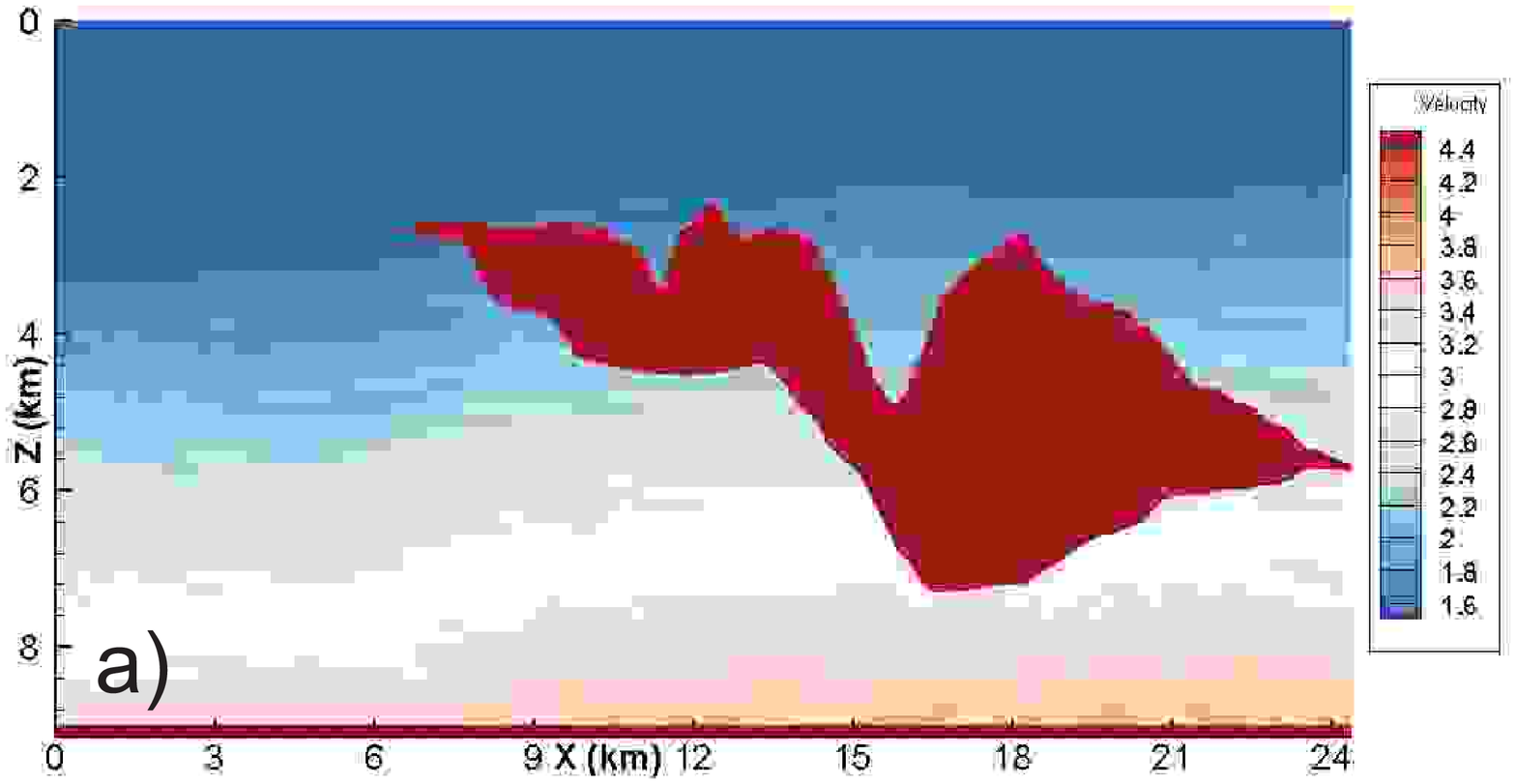}
        \end{subfigure}
              ~
                \begin{subfigure}[b]{0.49\textwidth}
                \includegraphics[width=\textwidth]{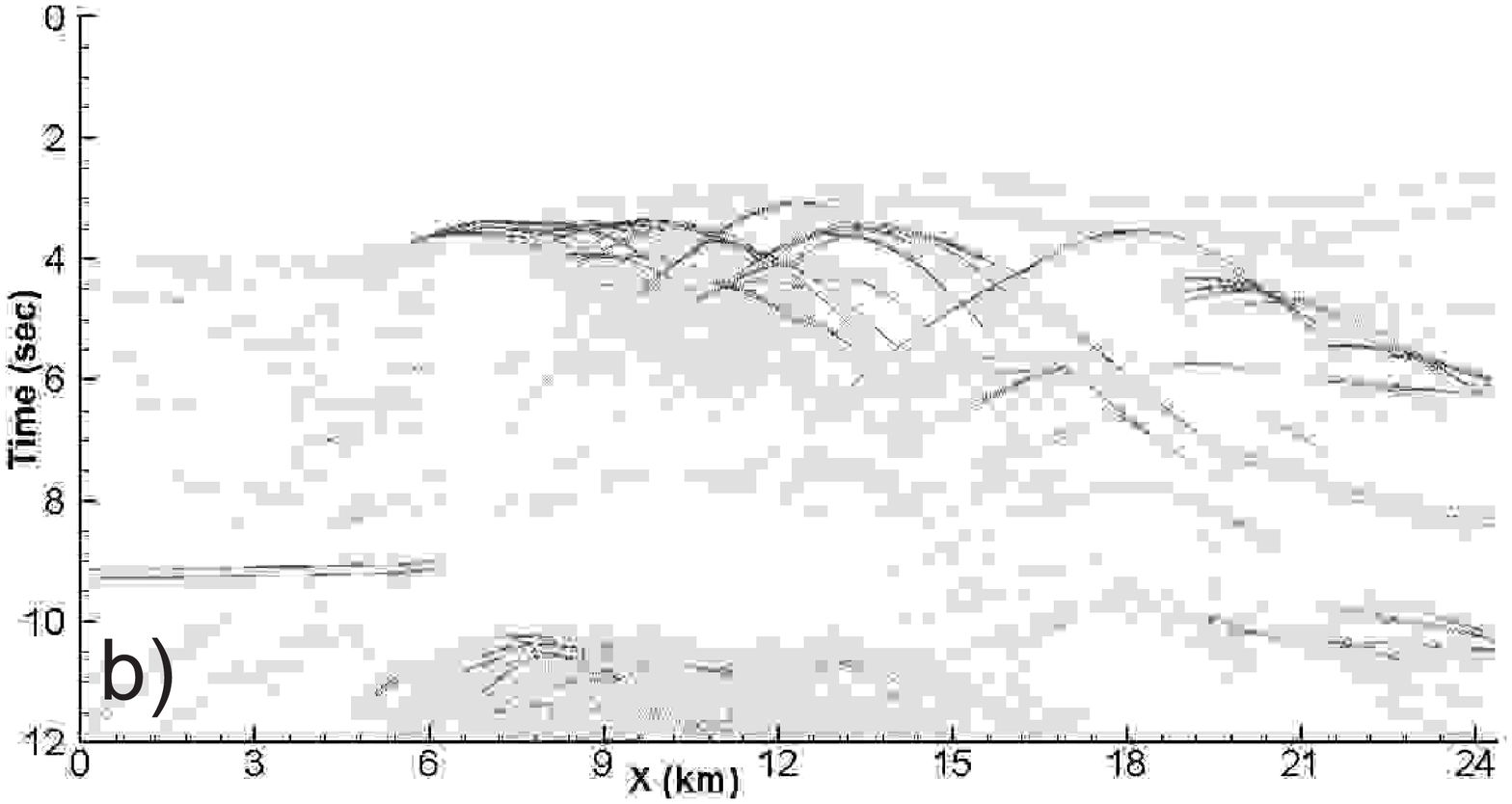}
        \end{subfigure}%
        \caption{(a) The 2D Sigsbee2A salt model and (b) zero-offset section.}
                \label{fig:sigsbee}
\end{figure}

For Sigsbee2Amodel \cite{Paffenholz} (Fig.~\ref{fig:sigsbee}a), theoretical seismograms (Fig.~\ref{fig:sigsbee}b) were calculated using the algorithm of explosive boundaries implemented in the Madagascar package \cite{Fomel2013}. For setting the boundary condition, the function for zero-offset section$\left.u(x,z,t)\right|_{z=0}=g(x,t)$ was expanded in series (\ref{series_lag}) for $n=3500$, parameters $\alpha=0$ and $\eta=300$ for $t\in[0,12]\;s$. The calculations were done on the meshes with the steps $h_x=h_z=7.62\;m$, $h_{x}=7.62\;m$  and $\;h_{z}=3.81\;m$.

According to the results of computational experiments, the FD method appeared to be unstable for the given velocity model and selected numerical parameters. Different degrees of smoothing the velocity medium model did not allow attaining the calculation stability for this method, therefore the FD method within this test is excluded from further consideration.

\begin{figure}[!h]
\centering

        \begin{subfigure}[b]{0.49\textwidth}
                \includegraphics[width=\textwidth]{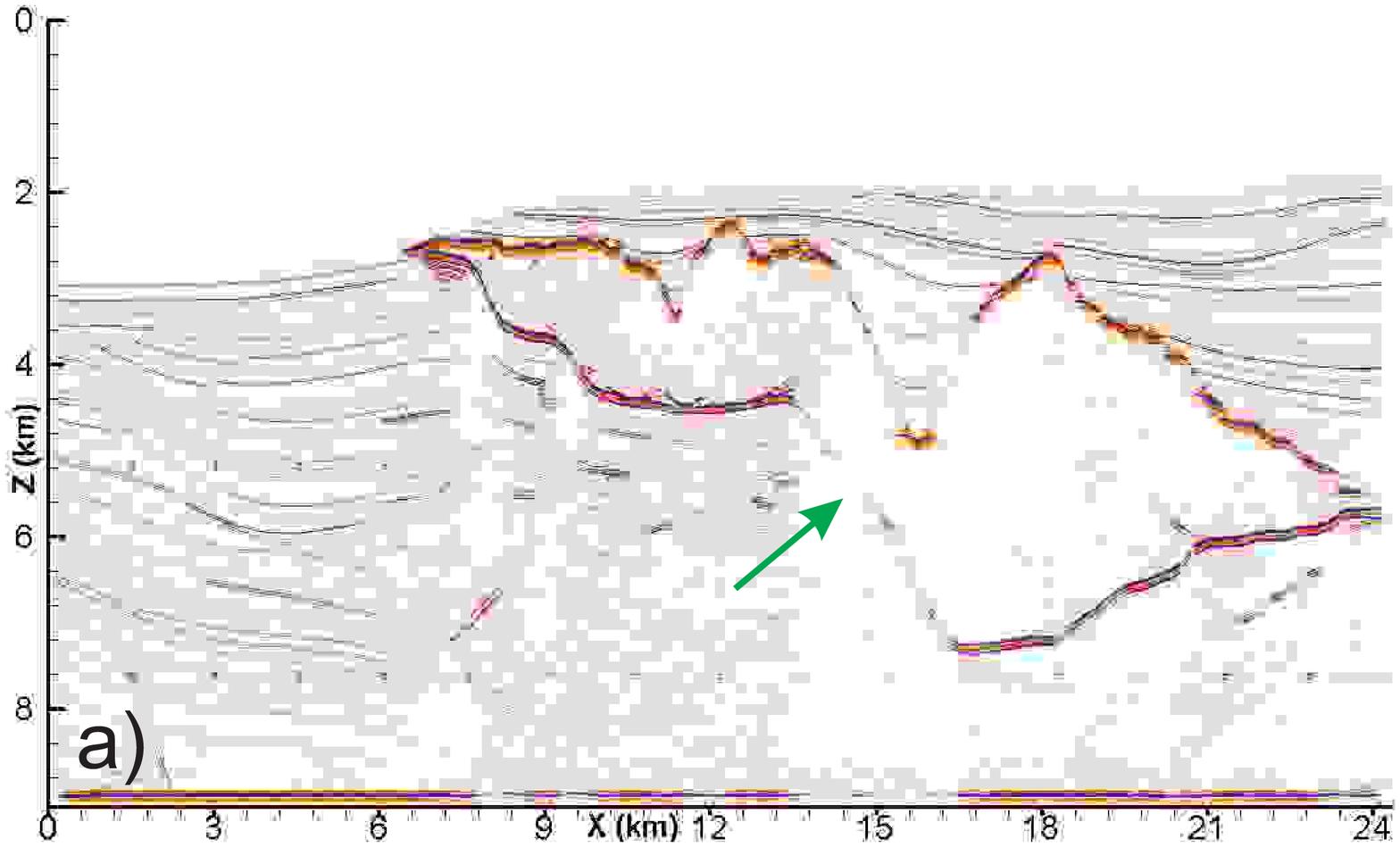}
                \label{fig:homo21}
        \end{subfigure}
        ~
        \begin{subfigure}[b]{0.49\textwidth}
                \includegraphics[width=\textwidth]{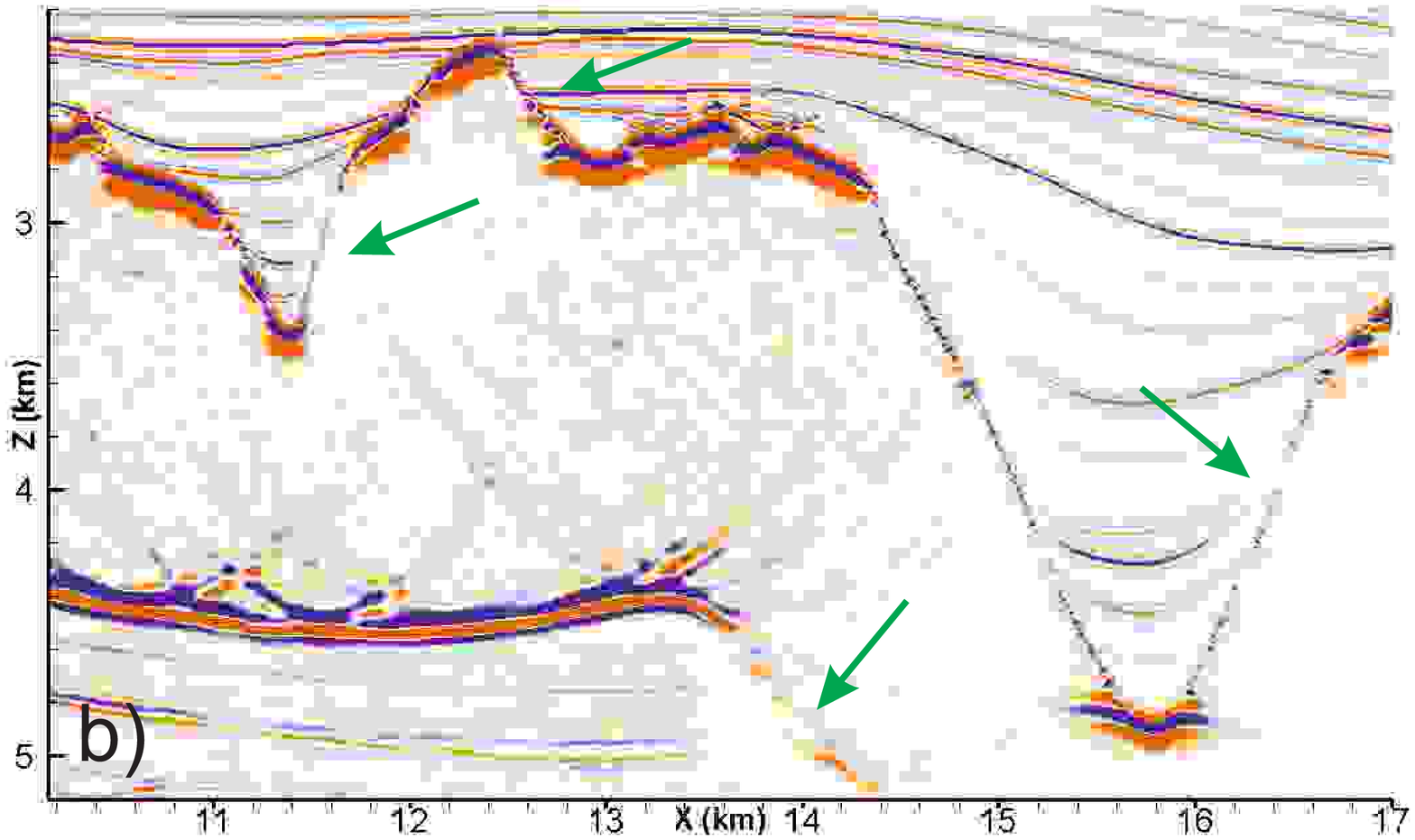}
                \label{fig:homo21}
        \end{subfigure}
        \caption{Snapshots for the wave field at $t=12\;s$ for the model and the zero-offset section in Fig.~\ref{fig:sigsbee} for the PSPI method for $h_x=h_z=7.62\;m$. Computational artifacts are indicated by arrows.}
         \label{fig:sig_pspi}
\end{figure}

\begin{figure}[!h]
\centering

        \begin{subfigure}[b]{0.49\textwidth}
                \includegraphics[width=\textwidth]{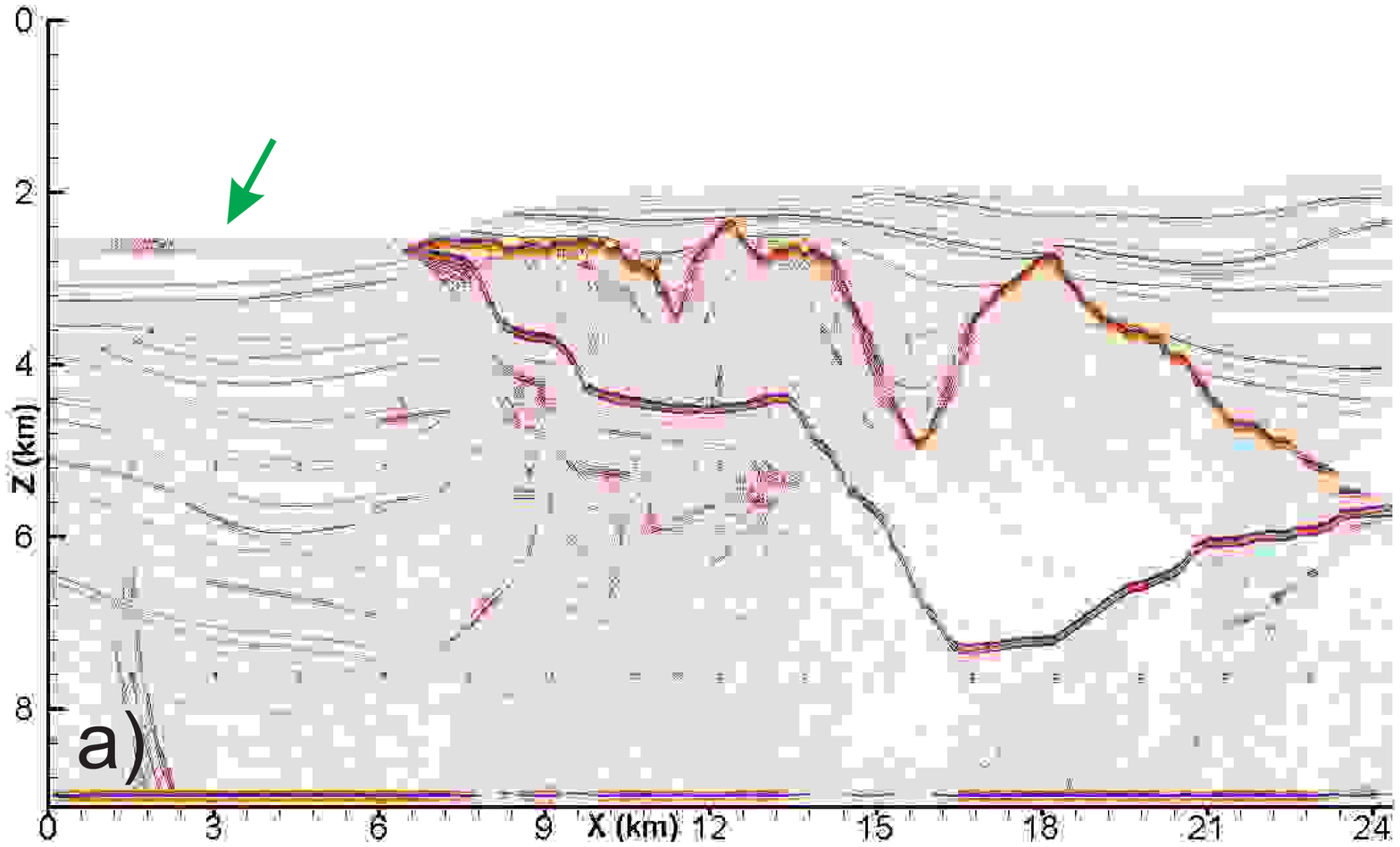}
                \label{fig:homo21}
        \end{subfigure}
          ~
                \begin{subfigure}[b]{0.49\textwidth}
                \includegraphics[width=\textwidth]{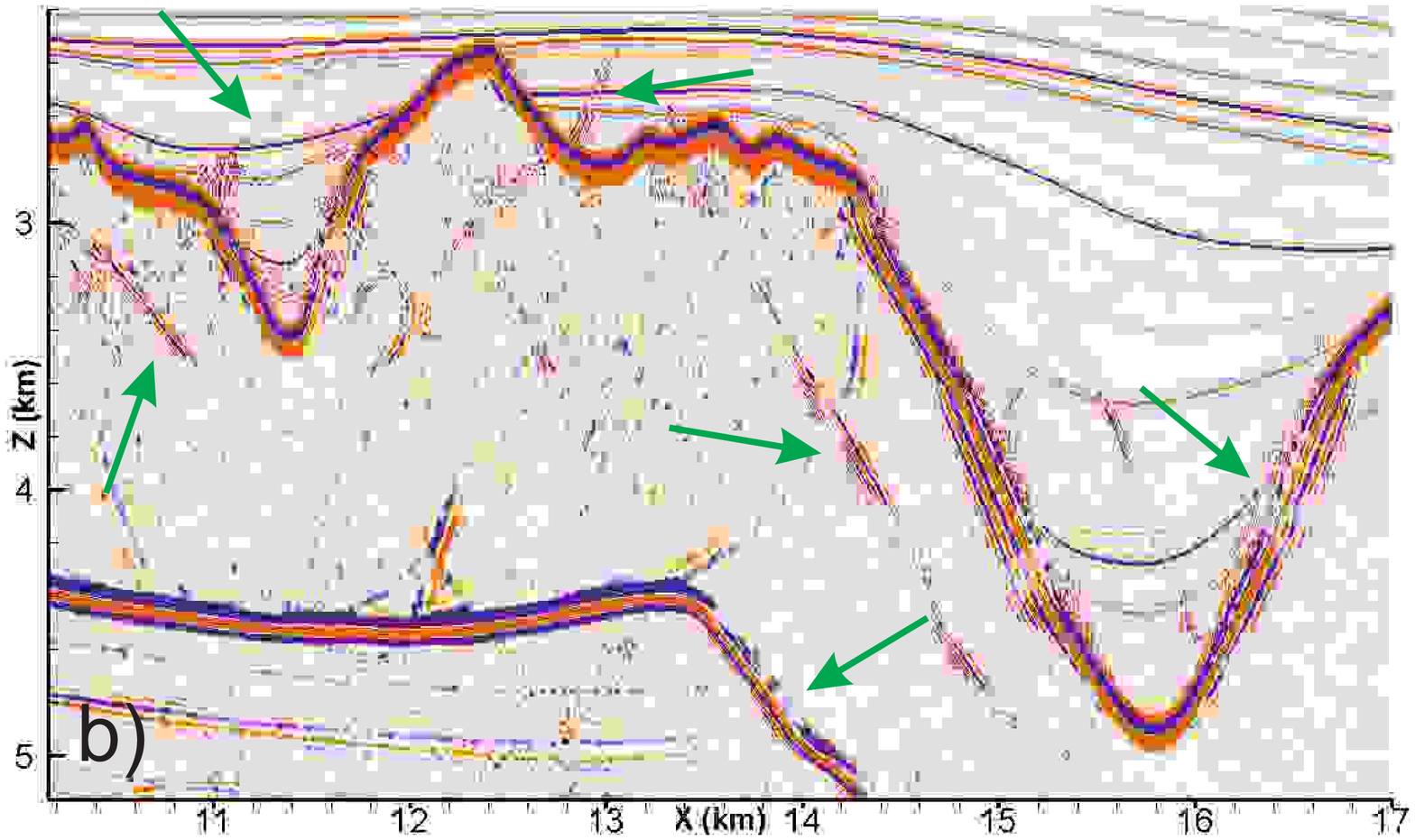}
                \label{fig:homo22}
        \end{subfigure}%
        \caption{Snapshots for the wave field at $t=12\;s$ for the model and the zero-offset section in Fig.~\ref{fig:sigsbee} for the FFD method for $h_x=h_z=7.62\;m$.  Computational artifacts are indicated by arrows.}
        \label{fig:sig_ffd}
\end{figure}

The PSPI method from the Seismic Unix package has shown unsatisfactory results because of bad energy focusing (Fig.~\ref{fig:sig_pspi}a). This method was unable to correctly map boundaries of the complex salt shape (Fig.~\ref{fig:sig_pspi}b). However note that the point type diffractors underlying the salt inclusion are correctly focused. The FFD method as compared to the PSPI has allowed obtaining a clearer image of the salt inclusion (Fig.~\ref{fig:sig_ffd}b), but the image additionally contains numerous artifacts (Fig.~\ref{fig:sig_ffd}a) due to approximate character of calculation formulas of the FFD method. The fact is, as well as in the previous test, an increase in resolving power of the mesh does not bring about an increase in calculation accuracy.

\begin{figure}[!htb]
\centering

        \begin{subfigure}[b]{0.49\textwidth}
                \includegraphics[width=\textwidth]{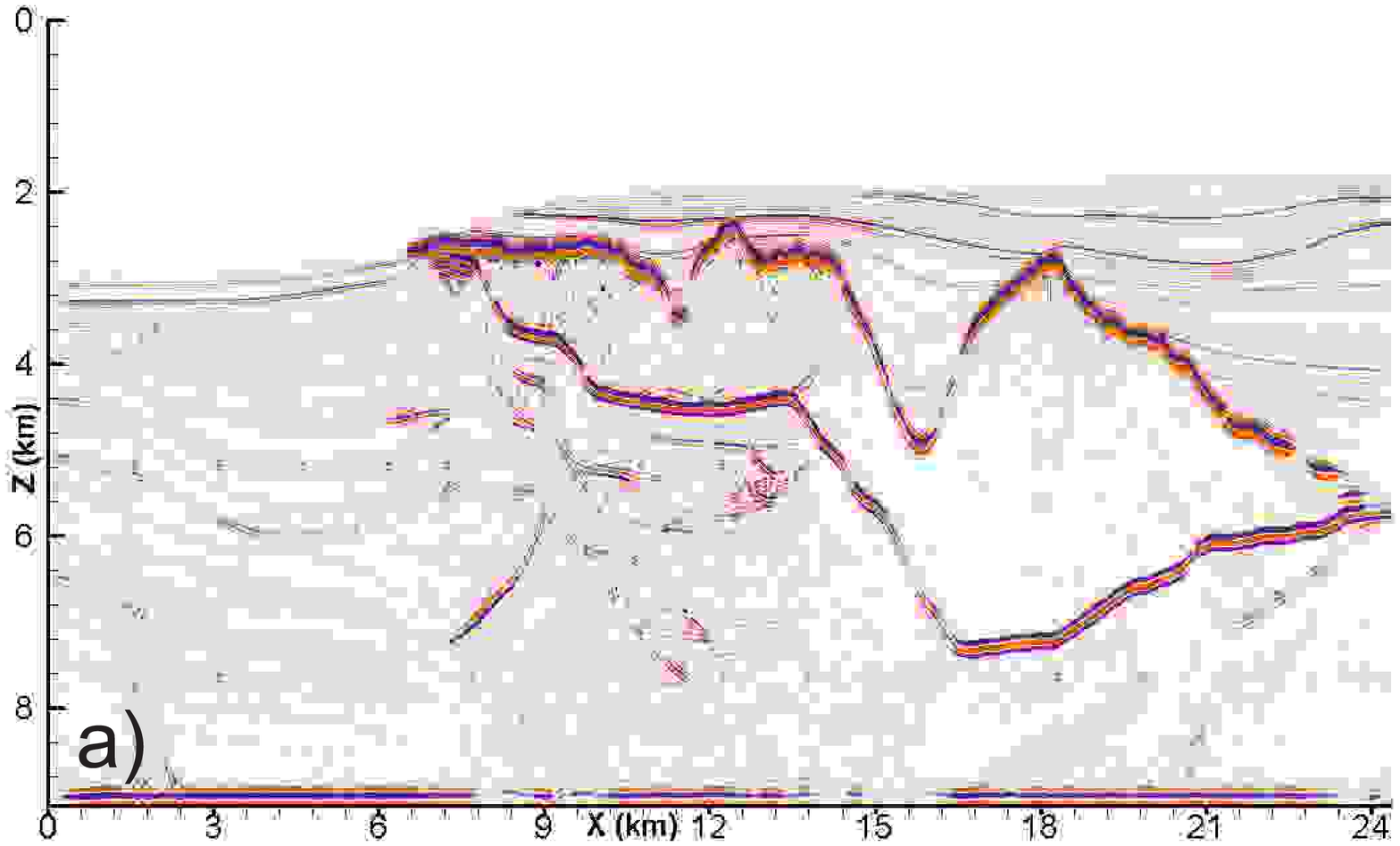}
                \label{fig:homo21}
        \end{subfigure}
             ~
                \begin{subfigure}[b]{0.49\textwidth}
                \includegraphics[width=\textwidth]{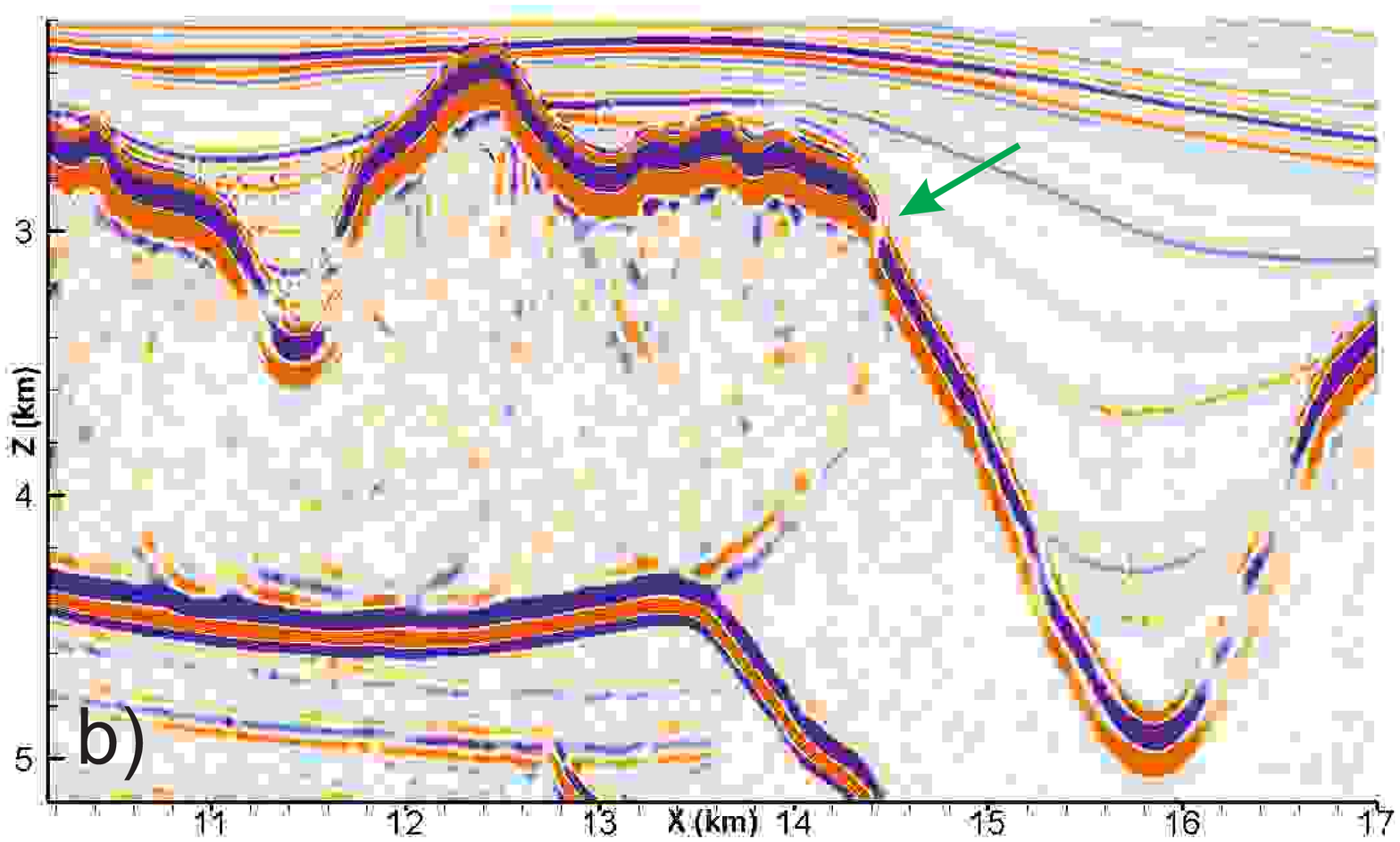}
                \label{fig:homo22}
        \end{subfigure}%

        \begin{subfigure}[b]{0.49\textwidth}
                \includegraphics[width=\textwidth]{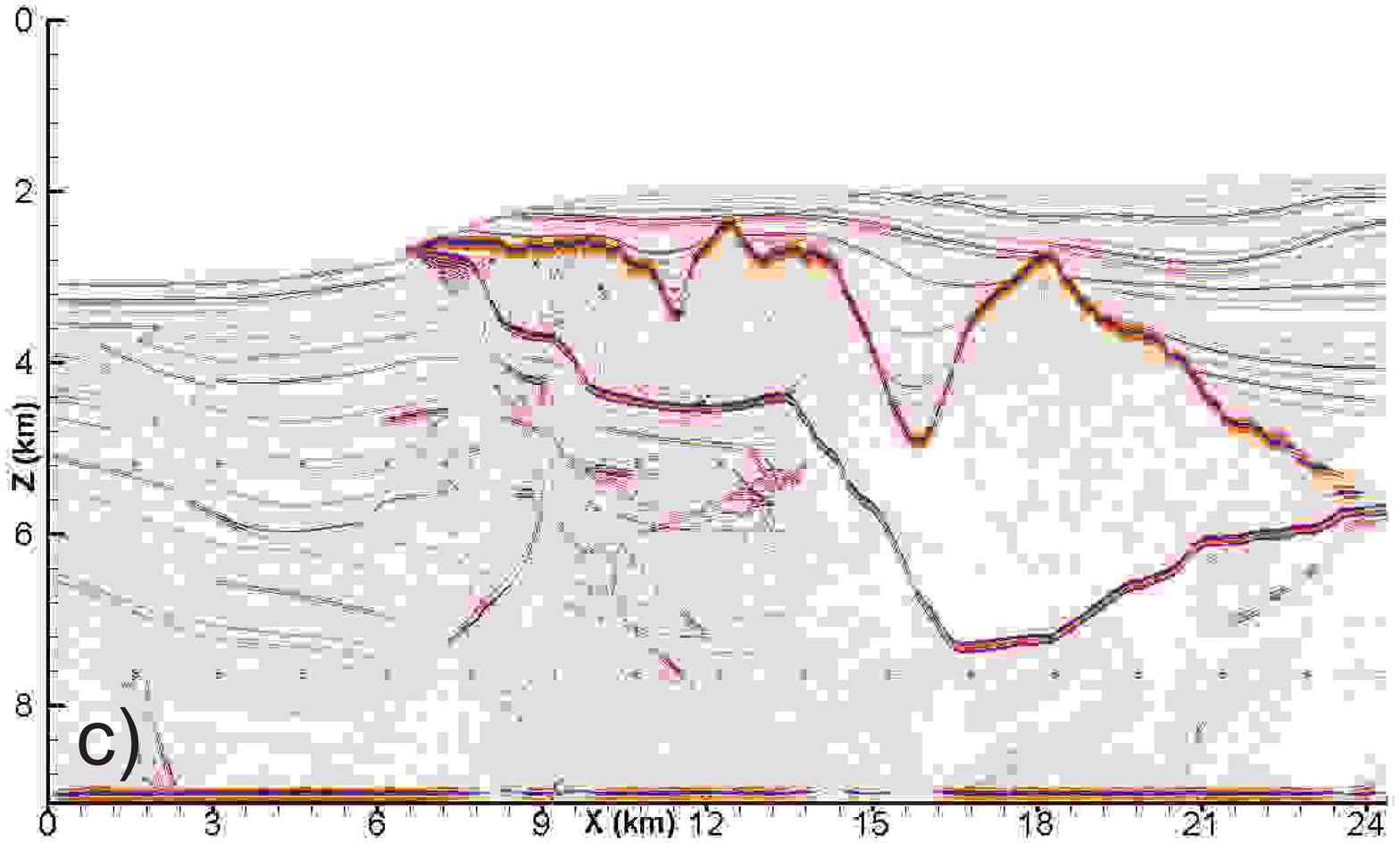}
                \label{fig:homo21}
        \end{subfigure}
        ~
        \begin{subfigure}[b]{0.49\textwidth}
                \includegraphics[width=\textwidth]{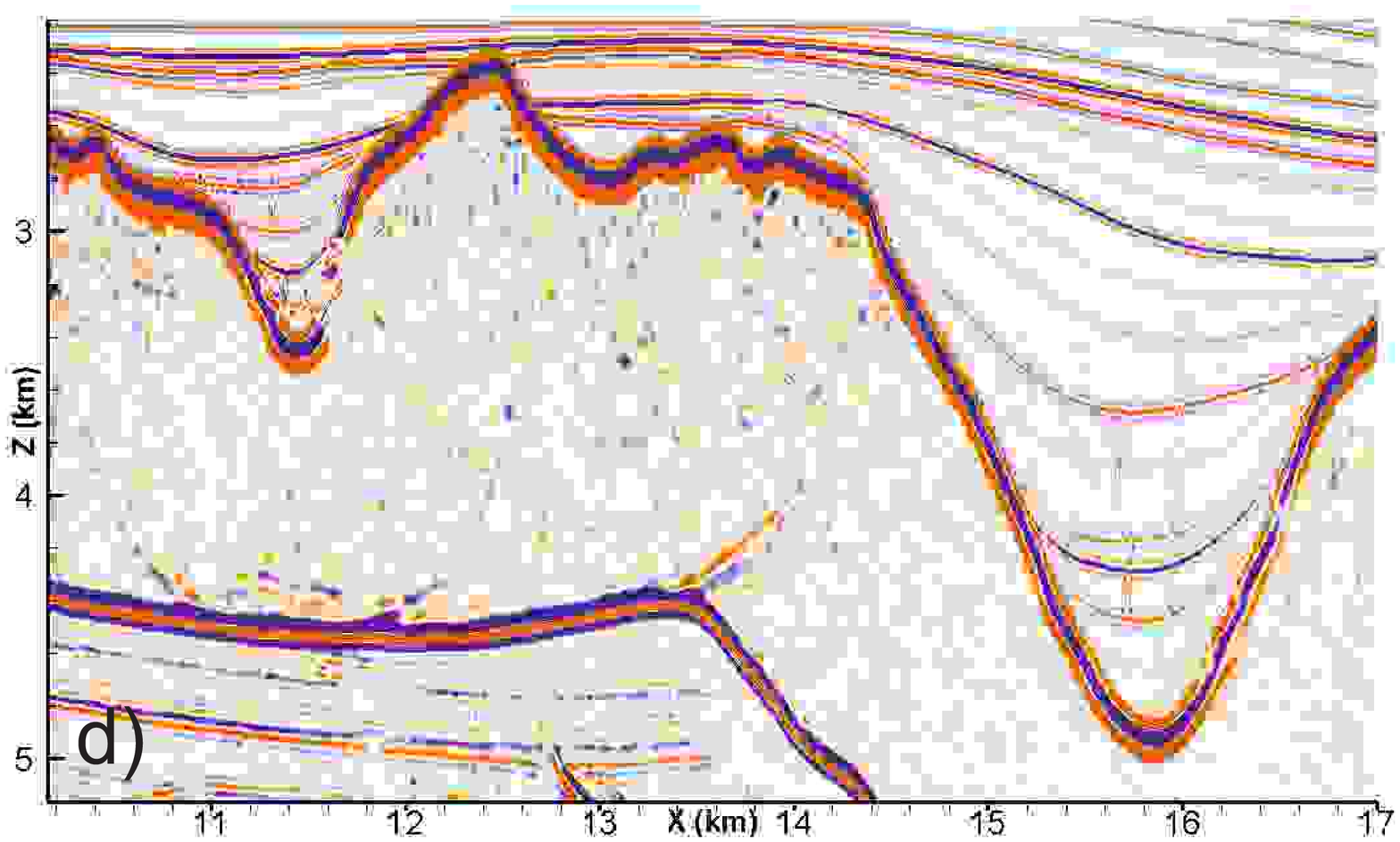}
                \label{fig:homo21}
        \end{subfigure}
\caption{Snapshots for the wave field at $t=12\;s$ for the model and the zero-offset section in Fig.~\ref{fig:sigsbee} for the LFD method (a),(b)~$h_x=h_z=7.62\;m$ and (c),(d)~$h_x=7.62,h_z=3.81\;m$. Computational artifacts are indicated by arrows.}
                \label{fig:sig_lfd}
\end{figure}

As opposed to the FFD and the PSPI methods, a conceptually different situation arises with the use of the LFD method: as is seen in Fig.~\ref{fig:sig_lfd}, the LFD method does not contain additional artifacts but only those which are simultaneously present for the FFD and the PSPI algorithms, i.e. the artifacts being present in input data of the problem (zero-offset section). Using the DRP approach to determine coefficients of the finite difference approximation $\partial^2/\partial x^2$  makes possible to attain high accuracy of calculation and to decrease numerical dispersion for $x$-direction. For $z$-direction, typical is numerical dissipation (Fig.~\ref{fig:sig_lfd}a,b) that brings about smoothing an image. However, the numerical dissipation can be decreased at the cost of choosing a smaller step $h_z$ (Fig.~\ref{fig:sig_lfd}c,d). Fourth approximation order due to the Richardson extrapolation allows the use of admissible steps for attaining the required calculation accuracy. Thus, errors caused by the LFD algorithm are controlled by choosing the step of a spatial mesh, whereas the choice of computational parameters for the PSPI and the FFD methods is a more ambiguous task, and decreasing a spatial mesh step for these methods does not really affect the imagery quality.

Thus, in the context of the post-stack migration using the method in question allows obtaining a more accurate image of the Earth's interior as compared to widespread mesh algorithms of solving the one-way equation. The method proposed can be also implemented for the pre-stack migration with a view to obtain high quality images. As opposed to numerical-analytical approaches requiring fulfillment of Fourier transform both with respect to time and with respect to space, the use of purely difference spatial approximation will enable us in the future to carry out calculations, first, for non-uniform meshes and, second, to implement a curvilinear boundary for the upper surface. The convergence of the new algorithm of the mesh step makes possible to evaluate the solution accuracy on a sequence of imbedded meshes thus separating approximation errors from errors of setting input data of the problem. For most of numerical migration methods it is necessary to preliminarily smooth a velocity medium model in order to provide the calculation stability. This brings about errors in the wave field kinematics and presents a severe problem in choosing a smoothing procedure. As the computational experiments show, the proposed LFD method for considered models and chosen computational parameters does not require preliminary smoothing of the medium model function, which counts in favor of a higher degree of stability as compared to existing mesh methods.
\subsection{Performance consideration}
As a matter of fact, a higher accuracy is not attained without supplementary computer costs. The single-processor implementation of the method proposed has shown that the LFD algorithm requires from $5$ to $10$ times more computational time than the PSPI, the FD, the FFD algorithms. The costs of the LFD algorithm for solving the SLAEs with banded matrices are essential, but the possibility of using optimized library functions for solving the SLAEs allows an increase in the method performance on the whole. Higher computer costs of the LFD algorithm can be justified because the FD method does not provide any accuracy and is weakly stable, and the FFD and the PSPI methods solve a certain approximate problem but not the one-way equation thus generating multiple noise and artifacts of the wave field.
\begin{figure}[!h]
\centering
                \includegraphics[width=\textwidth]{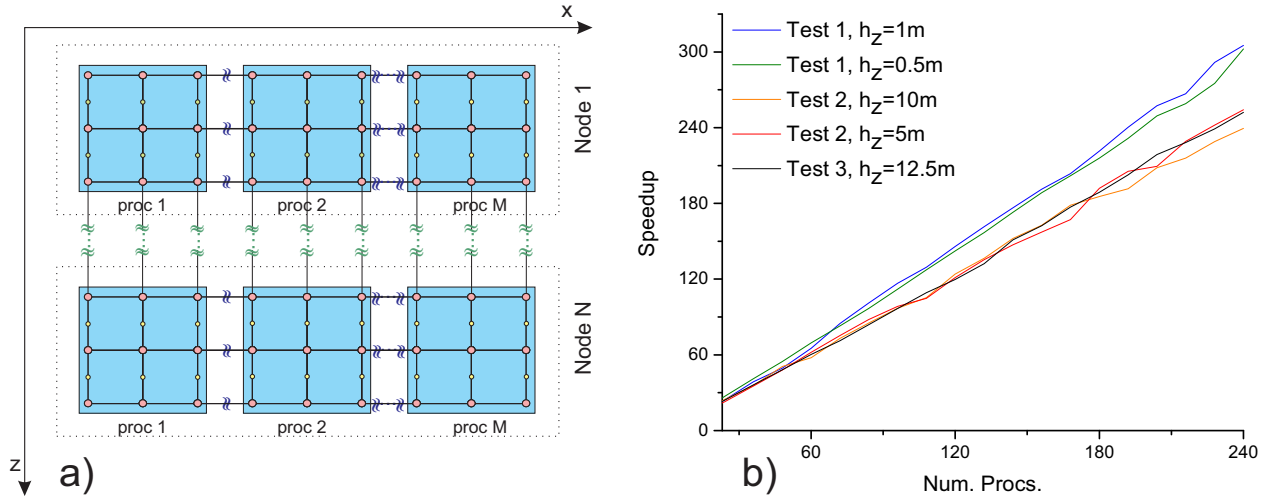}
\caption{(a)~Data decomposition among processors. (b)~Speedup versus the number of processors for different meshes and velocity models.}
\label{fig:parallel}

\end{figure}

For testing the performance of the method proposed the calculation time for a different number of mesh nodes and the number of processors was assessed. The use of multi-processor systems makes absolute time costs insignificant because with a maximum number of processors the calculation time for all the test problems does not exceed several seconds (Table~\ref{table33}).
\begin{table}[!h]
\center \small
\begin{tabular}{cccccccccc}
  \hline
     & \multicolumn{1}{c}{$Test\;1$}& \multicolumn{1}{c}{$ Test\;1$}& \multicolumn{1}{c}{$Test\;2$}& \multicolumn{1}{c}{$Test\;2$}& \multicolumn{1}{c}{$Test\;3$}\\
     $N_x \times N_z$ & \multicolumn{1}{c}{$1624\times1024$}& \multicolumn{1}{c}{$3248\times2048$}& \multicolumn{1}{c}{$462\times300$}& \multicolumn{1}{c}{$924\times600$}& \multicolumn{1}{c}{$6444\times4800$}\\
    num proc.
    \\ \hline

  12  &247&1008&44 &119&944  \\
  24  &122&446&21&62&472 \\
  48 &76&294&14.8&41&313 \\
  ...  &...&...&... &...&...  \\
  216  &11&46&2.4 &6.2&49.6  \\
  228  &10&43&2.3&5.9&47.3 \\
  240 &9.7&40&2.2&5.6&44.9 \\
  \hline
  \end{tabular}
\caption{Calculation time (seconds) versus the number of processors for meshes of different resolution for all the tests.}
\label{table33}
\end{table}
The dependence of the speedup value on the number of processors for all the tests is presented in Fig.~\ref{fig:parallel}b, indicating to the linear character of the dependence, for the first test attaining the super-linear speedup. A similar effect of super-linear speedup was observed when using the parallel dichotomy algorithm within the implementation of the alternating implicit direction method for solving the 2D Poisson equation\cite{terekhov:Dichotomy}. Thus, the dichotomy algorithm implemented in the context of the approach proposed for solving the 2D one-way equation allows the effective use of a considerable number of processors. We can conclude that the chosen data decomposition of the problem (Fig.~\ref{fig:parallel}b) is fairly appropriate, as multiple exchanges for x-direction are instantly made at the expense of using the shared memory within one computational node, while slow inter-node communicative interactions for $z$-direction for solving problem (\ref{reduce_eq}) are done only once.
\section{Conclusion}
One of the major results of this research is in  managing to construct a stable, efficient, purely finite difference method of high accuracy order for solving the one-way equation. By now this has presented considerable difficulties in numerical analysis due to the necessity of providing high accuracy and at the same time the calculation stability. Replacing Fourier transform by the Laguerre transform  after employing the finite difference approximation of spatial derivatives, using algebraic transformations made possible to obtain the banded SLAEs appropriate for solutions by direct methods, in particular, by the parallel dichotomy algorithm. The rejection from using the splitting of the Marchuk-Strang type for the problem operator has excluded appearance of numerous computing artifacts and has increased both the order of approximation and the stability of the method as a whole. The one-way equation does not allow a correct description of the wave field amplitude for inhomogeneous media, hence the optimization of values of finite difference scheme coefficients aimed at decreasing the phase errors is a reasonable approach. The numerical experiments have shown that such DRP schemes make possible to decrease the mesh step in horizontal approximately by the factor of two as compared to classical difference schemes based on the Taylor series. Implementation of the Richardson extrapolation procedure has not only increased the accuracy of the algorithm up to fourth order relative to the mesh step in depth, but also has provided the stability of the method as compared to using schemes of the Adams-Moulton type. Applying implicit Runge-Kutta methods of high orders of accuracy is unreasonable, as in this case it is required to solve SLAEs of much higher dimensions. In addition, as well as for the Richardson extrapolation, the interpolation of the wave field values for auxiliary mesh nodes will  be needed.

The algorithm proposed can be generated to solving 3D problems as well. To this end, the efficient iterative procedure for solving SLAEs should be developed. This is because the parallel dichotomy algorithm, being a direct method, is not intended for banded matrices of a large width. Also, it is necessary to reject not only the splitting of the Marchuk-Strang type for the problem operator, but also the splitting along horizontal directions that is used for decreasing computer costs, but generates numerical anisotropy.  The efficient algorithm for solving the 3D one-way equation will be presented in the next publication in the near future.
\section{Acknowledgments.}
The author is grateful to Gennady Zhernyak for fruitful discussions of the paper.
\newpage
\bibliography{base}

\begin{thebibliography}{10}

\bibitem{Claerbout1971}
J.F. Claerbout.
\newblock Toward a unified theory of reflector mapping.
\newblock {\em Geophysics}, 36(3):467--481, 1971.

\bibitem{seis_interpr}
O.~Yilmaz and S.M. Doherty.
\newblock {\em Seismic Data Analysis: Processing, Inversion, and Interpretation
  of Seismic Data}.
\newblock Society of Exploration Geophysicist, 2001.

\bibitem{Angus2013}
D.A. Angus.
\newblock The one-way wave equation: A full-waveform tool for modeling seismic
  body wave phenomena.
\newblock {\em Surveys in Geophysics}, 35(2):359--393, 2013.

\bibitem{Tappert1977}
F.D. Tappert.
\newblock {\em Wave Propagation and Underwater Acoustics}, chapter The
  parabolic approximation method, pages 224--287.
\newblock Springer Berlin Heidelberg, Berlin, Heidelberg, 1977.

\bibitem{Jensen2011}
F.B. Jensen, W.A. Kuperman, M.B. Porter, and H.~Schmidt.
\newblock {\em Computational Ocean Acoustics}.
\newblock Springer-Verlag New York, 2011.

\bibitem{Lee2000}
D.~Lee, A.D. Pierce, and Er-C. Shang.
\newblock Parabolic equation development in the twentieth century.
\newblock {\em Journal of Computational Acoustics}, 08(04):527--637, 2000.

\bibitem{Lindman1975}
E.L. Lindman.
\newblock ''{F}ree-space'' boundary vonditions for the time dependent wave
  equation.
\newblock {\em Journal of Computational Physics}, 18(1):66 -- 78, 1975.

\bibitem{Engquist1977}
B.~Engquist and A.~Majda.
\newblock Absorbing boundary conditions for the numerical simulation of waves.
\newblock {\em Math. Comp.}, 31:629--651, 1977.

\bibitem{Trefethen1986}
L.N. Trefethen and L.~Halpern.
\newblock Well-posedness of one-way wave equations and absorbing boundary
  conditions.
\newblock {\em Mathematics of Computation}, 47(176):421--435, 1986.

\bibitem{Leontovich1946}
M.~Leontovich and V.~Fock.
\newblock Solution of the problem of propagation of electromagnetic waves along
  the earth’s surface by the method of parabolic equation.
\newblock {\em Acad. Sci. USSR J. Phys.}, 10:13--24, 1946.
\newblock (Engl. transl).

\bibitem{Claerbout:1985}
J.F. Claerbout.
\newblock {\em Imaging the Earth's Interior}.
\newblock Blackwell Scientific Publications, Inc., Cambridge, MA, USA, 1985.

\bibitem{Lee1982}
D.~Lee and G.~Botseas.
\newblock Ifd: An implicit finite-difference computer model for solving the
  parabolic equation.
\newblock Technical Report 6659, Nav. Underwater Syst. Ctr, 1982.

\bibitem{Halpern1988}
L.~Halpern and L.N. Trefethen.
\newblock Wide-angle one-way wave equations.
\newblock {\em J. Acoust. Soc. Am.}, 84(4):1397--404, 1988.

\bibitem{Bamberger1988}
A.~Bamberger, B.~Engquist, L.~Halpern, and P.~Joly.
\newblock Higher order paraxial wave equation approximations in heterogeneous
  media.
\newblock {\em SIAM Journal on Applied Mathematics}, 48(1):129--154, 1988.

\bibitem{Bamberger1988a}
A.~Bamberger, B.~Engquist, L.~Halpern, and P.~Joly.
\newblock Parabolic wave equation approximations in heterogenous media.
\newblock {\em SIAM Journal on Applied Mathematics}, 48(1):99--128, 1988.

\bibitem{Marchuk1968}
G.I. Marchuk.
\newblock Some application of splitting-up methods to the solution of
  mathematical physics problems.
\newblock {\em Applik Mat}, 13(2):103--132, 1968.

\bibitem{Strang1968}
G.~Strang.
\newblock On the construction and comparison of difference schemes.
\newblock {\em SIAM Journal on Numerical Analysis}, 5(3):506--517, 1968.

\bibitem{Marchuk1990}
G.I. Marchuk.
\newblock Handbook of numerical analysis.
\newblock In P.~G. Ciarlet and J.~L. Lions, editors, {\em Handbook of Numerical
  Analysis}, volume~1, pages 197--460. Elsevior Science Publishers BV, 1990.

\bibitem{Bunks1995}
C.~Bunks.
\newblock Effective filtering of artifacts for implicit finite-difference
  paraxial wave equation migration1.
\newblock {\em Geophysical Prospecting}, 43(2):203--220, 1995.

\bibitem{Fei1995}
T.~Fei and K.~Larner.
\newblock Elimination of numerical dispersion in finite-difference modeling and
  migration by flux-corrected transport.
\newblock {\em Geophysics}, 60(6):1830--1842, 1995.

\bibitem{Gazdag1978}
J.~Gazdag.
\newblock Wave equation migration with the phase-shift method.
\newblock {\em Geophysics}, 43(7):1342--1351, 1978.

\bibitem{Gazdag1984}
J.~Gazdag and P.~Sguazzero.
\newblock Migration of seismic data by phase shift plus interpolation.
\newblock {\em Geophysics}, 49(2):124--131, 1984.

\bibitem{Ristow01121994}
D.~Ristow and T.~R\"{u}hl.
\newblock Fourier finite-difference migration.
\newblock {\em Geophysics}, 59(12):1882--1893, 1994.

\bibitem{Guddati2006}
M.N. Guddati.
\newblock Arbitrarily wide-angle wave equations for complex media.
\newblock {\em Computer Methods in Applied Mechanics and Engineering}, 195:65
  -- 93, 2006.

\bibitem{abramowitz+stegun}
M.~Abramowitz and I.A. Stegun.
\newblock {\em Handbook of Mathematical Functions with Formulas, Graphs, and
  Mathematical Tables}.
\newblock Dover, New York, ninth dover printing, tenth gpo printing edition,
  1964.

\bibitem{Mikh2003}
B.G. Mikhailenko.
\newblock Simulation of seismic wave propagation in heterogeneous media.
\newblock {\em Siberian J. of Numer. Mathematics}, 6:415--429, 2003.

\bibitem{fatab2011}
A.G. Fatyanov and A.V. Terekhov.
\newblock High-performance modeling acoustic and elastic waves using the
  parallel dichotomy algorithm.
\newblock {\em J. Comp. Phys.}, 230(5):1992--2003, 2011.

\bibitem{Terekhov:2013}
A.V. Terekhov.
\newblock A fast parallel algorithm for solving block-tridiagonal systems of
  linear equations including the domain decomposition method.
\newblock {\em Parallel Comput.}, 39(6-7):245--258, 2013.

\bibitem{Terekhov2015206}
A.V. Terekhov.
\newblock Spectral-difference parallel algorithm for the seismic forward
  modeling in the presence of complex topography.
\newblock {\em Journal of Applied Geophysics}, 115(0):206--219, 2015.

\bibitem{Lee01101985}
Myung~W. Lee and Sang~Y. Suh.
\newblock Optimization of one-way wave equations.
\newblock {\em Geophysics}, 50(10):1634--1637, 1985.

\bibitem{Integral_Transform}
L.~Debnath and D.~Bhatta.
\newblock {\em Integral Transforms and Their Applications, Second Edition}.
\newblock Taylor \& Francis, 2006.

\bibitem{Mikhailenko1999}
B.G. Mikhailenko.
\newblock Spectral laguerre method for the approximate solution of time
  dependent problems.
\newblock {\em Applied Mathematics Letters}, 12:105--110, 1999.

\bibitem{CrankNicolson2}
J.~Crank and P.~Nicolson.
\newblock A practical method for numerical evaluation of solutions of partial
  differential equations of the heat-conduction type.
\newblock {\em Proc. Camb. Phil. Soc.}, 43(1):50--67, 1947.

\bibitem{Zhang2005}
Y.~Zhang, G.~Zhang, and N.~Bleistein.
\newblock Theory of true-amplitude one-way wave equations and true-amplitude
  common-shot migration.
\newblock {\em Geophysics}, 70(4):E1--E10, 2005.

\bibitem{Tam1993262}
C.K.W. Tam and J.C. Webb.
\newblock Dispersion-relation-preserving finite difference schemes for
  computational acoustics.
\newblock {\em Journal of Computational Physics}, 107(2):262 -- 281, 1993.

\bibitem{GPR:GPR12160}
W.~Liang, Y.~Wang, and C.~Yang.
\newblock Determining finite difference weights for the acoustic wave equation
  by a new dispersion-relationship-preserving method.
\newblock {\em Geophysical Prospecting}, 63(1):11--22, 2015.

\bibitem{Chu2012}
C.~Chu and P.L. Stoffa.
\newblock Determination of finite-difference weights using scaled binomial
  windows.
\newblock {\em Geophysics}, 77(3):W17--W26, 2012.

\bibitem{Zhang2013511}
J.-H. Zhang and Z.-X. Yao.
\newblock Optimized explicit finite-difference schemes for spatial derivatives
  using maximum norm.
\newblock {\em Journal of Computational Physics}, 250:511 -- 526, 2013.

\bibitem{Samarskii2001}
A.A. Samarskii.
\newblock {\em The Theory of Difference Schemes}.
\newblock Marcel Dekker, 2001.

\bibitem{Cottle1974}
R.W. Cottle.
\newblock Manifestations of the schur complement.
\newblock {\em Linear Algebra and its Applications}, 8(3):189 -- 211, 1974.

\bibitem{Henderson1981}
H.V. Henderson and S.R. Searle.
\newblock On deriving the inverse of a sum of matrices.
\newblock {\em SIAM Review}, 23(1):53--60, 1981.

\bibitem{rohtua}
L.F. Richardson.
\newblock On the approximate arithmetical solution by finite differences of
  physical problems involving differential equations, with an application to
  the stresses in a masonry dam.
\newblock {\em Proceedings of the Royal Society of London A: Mathematical,
  Physical and Engineering Sciences}, 83(563):335--336, 1910.

\bibitem{Marchuk:99314}
G.~I. Marchuk and V.~V. Shaidurov.
\newblock {\em Difference methods and their extrapolations}.
\newblock Appl. Math. Springer, New York, NY, 1983.

\bibitem{terekhov:Dichotomy}
A.V. Terekhov.
\newblock Parallel dichotomy algorithm for solving tridiagonal system of linear
  equations with multiple right-hand sides.
\newblock {\em Parallel Comput.}, 36(8):423--438, 2010.

\bibitem{Terekhov2}
A.V. Terekhov.
\newblock A highly scalable parallel algorithm for solving toeplitz tridiagonal
  systems of linear equations.
\newblock {\em Journal of Parallel and Distributed Computing}, 87:102--108,
  2016.

\bibitem{Gray2001}
S.H. Gray, J.~Etgen, J.~Dellinger, and D.~Whitmore.
\newblock Seismic migration problems and solutions.
\newblock {\em Geophysics}, 66(5):1622--1640, 2001.

\bibitem{Biondi2006}
B.~Biondi.
\newblock {\em 3D Seismic Imaging}.
\newblock Society of Exploration Geophysicists, 2006.

\bibitem{GaussianBeams_popov}
M.M. Popov.
\newblock A new method of computation of wave fields using gaussian beams.
\newblock {\em Wave Motion}, 4, 1982.

\bibitem{Cerveny1985}
V.~Cerveny.
\newblock Gaussian beam synthetic seismograms.
\newblock {\em J. Geophys.}, 58:44--72, 1985.

\bibitem{Paffenholz}
J.~Paffenholz, B.~McLain, J.~Zaske, and P.J. Keliher.
\newblock {\em Subsalt multiple attenuation and imaging: Observations from the
  Sigsbee2B synthetic dataset}, chapter 538, pages 2122--2125.

\bibitem{Fomel2013}
S.~Fomel, P.~Sava, I.~Vlad, Y.~Liu, and V.~Bashkardin.
\newblock Madagascar: open-source software project for multidimensional data
  analysis and reproducible computational experiments.
\newblock {\em Journal of Open Research Software 1:e8}, 2013.

\end{thebibliography}
\end{document}